\documentclass[11pt,letterpaper]{amsart}

\usepackage{amsmath,amssymb,amsthm,verbatim, csquotes,mathrsfs,stmaryrd}
\usepackage{hyperref}
\usepackage{xcolor}
\usepackage{esint}
\usepackage{mathtools}
\usepackage{graphicx}
\usepackage{float}
\usepackage{caption}
\usepackage{enumitem}
\mathtoolsset{showonlyrefs}

\usepackage[top=3.5cm, bottom=4cm, right=2.5cm, left=2.5cm]{geometry}

\allowdisplaybreaks
\mathtoolsset{showonlyrefs=true}

\usepackage{geometry}
\newtheorem{theorem}{Theorem}[section]
\newtheorem{proposition}[theorem]{Proposition}
\newtheorem{lemma}[theorem]{Lemma}
\newtheorem{corollary}[theorem]{Corollary}

\newtheorem{example}[theorem]{Example}
\theoremstyle{remark}
\newtheorem{remark}[theorem]{Remark}
\theoremstyle{definition}

\numberwithin{equation}{section}

\newcommand{\vol}{\operatorname{vol}}
\newcommand{\Ric}{\operatorname{Ric}}
\newcommand{\Ricz}{\mathring{\Ric}}
\newcommand{\dv}{\,dv}
\newcommand{\hesz}{\mathring{\nabla}^2}
\newcommand{\R}{\mathbb{R}}
\newcommand{\Z}{\mathbb{Z}}
\newcommand{\dvg}{\operatorname{div}}

\newcommand{\tr}{\operatorname{tr}}
\newcommand{\Sn}{\mathbb{S}^n}
\newcommand{\Sph}{ {\mathbb{S}}}

\title[Sharp integral Ricci pinching]
{Sharp integral Ricci pinching on manifolds\\ with positive total scalar curvature}
\author{Guofang Wang}
\address[G. Wang]{Mathematisches Institut, Albert-Ludwigs-Universit\"at Freiburg, Ernst-Zermelo-Str.~1, 79104 Freiburg, Germany}
\email{guofang.wang@math.uni-freiburg.de}

\author{Wei Wei}
\address[W.~Wei]{Nanjing University, School of Mathematics, Nanjing 210093, P.R. China}
\email{wei\_wei@nju.edu.cn}
\date{23 November 2026}

\subjclass[2020]{Primary 53C21; Secondary 53C23, 58J60}

\keywords{Ricci curvature, integral pinching, scalar curvature, sharp geometric inequality}

\begin{document}

\begin{abstract}
Let $(M,g)$ be a closed connected Riemannian three-manifold with $\int_M R\,\dv>0$, and assume that $M$ is not diffeomorphic to a spherical space form We prove the sharp integral Ricci pinching inequality
\begin{equation*}
\vol(g)\int_M|\Ricz|^2\dv
\ge \frac16\Bigl(\int_M R\dv\Bigr)^2.
\end{equation*}
Moreover, equality holds if and only if $(M,g)$ is locally isometric to $\mathbb S^2(a)\times\R$ for some $a>0$. 
By the same method, we prove $\vol(g)\int_M|\Ricz|^2\dv\ge\frac1{12}(\int_MR\dv)^2$ on closed four-manifolds with $\int_MR\dv>0$ that have a finite cover with positive first Betti number, and we provide new proof of Gursky's inequality.
We also show that the corresponding inequality fails in every dimension $n\ge5$, and we give a metric on $\mathbb S^2\times\mathbb S^1$ that disproves a conjecture of Ge, Li and Lian.

\end{abstract}

\maketitle

\section{Introduction}

Sharp inequalities relating global curvature quantities often characterize  model geometries as optimizers. An example is the classical Sobolev inequality on the round sphere $\Sph^n$, which states
\[
\int_{\Sn} R_g\,\dv_g \ge n(n-1)\vol(g_{\Sn})^{\frac2n}\,\vol(g)^{\frac{n-2}{n}},\qquad \forall\, g\in [g_{\Sn}],
\]
with equality if and only if $g$ is a round metric. The objective  of this paper is to obtain analogous optimal inequalities that characterize other model spaces, notably $\Sph^{n-1}\times\Sph^1$. Such characterizations are typically more subtle, since one needs additional hypotheses to rule out the standard sphere.

A very interesting result in this direction is due to Gursky~\cite{Gu}:

\smallskip
{\it If $M$ is a closed $4$-manifold with $b_1(M)>0$, then every metric $g$ with $R>0$ satisfies
\begin{equation}\label{eq:1.1}
    \int_M |\Ricz|^2 \ge \frac{1}{12}\int_M R^2,
\end{equation}
with equality if and only if $g$ is conformal to a metric locally isometric to $\Sph^{3}(a)\times \R$ for some $a>0$.}

\smallskip
In fact, it needs only the assumption that the Yamabe constant $Y(M,g)>0$ \cite{Gu}.
This result leads naturally to the following question in arbitrary dimension.

\smallskip
\noindent{\bf Question 1.}
{\it Let $M^n$ be a closed manifold with $b_1(M)>0$. Does every metric with $R>0$ satisfy
\begin{equation}\label{eq:1.2}
    \int_M |\Ricz|^2 \ge \frac{1}{n(n-1)}\int_M R^2,
\end{equation}
with equality if and only if $g$ is conformal to a metric locally isometric to $\Sph^{n-1}(a)\times \R$ for some $a>0$?
}

\smallskip

In dimension three, a related program has been carried out  to characterize $\Sph^2\times\Sph^1$ via integral inequalities, or integral pinching. Motivated by the $4$-dimensional conformal sphere theorem of Chang--Gursky--Yang~\cite{CGY03}, Ge--Lin--Wang~\cite{GLW} proved that any $3$-manifold admitting a metric with $R_g>0$ and $\int_M\sigma_2(g)\,\dv_g>0$ is a quotient of the $3$-sphere; the same conclusion was obtained independently by Catino--Djadli~\cite{CD}. 

The condition $\int_M\sigma_2>0$ is not conformal invariant in dimension $3$, and is not expected to be optimal. Here $\sigma_1(g)=\frac1{2(n-1)}R_g$ and $\sigma_2(g)=\frac{(\tr A_g)^2-|A_g|^2}{2}$, where $A_g=\frac1{n-2}(\Ric_g-\frac1{2(n-1)}R_g g)$,  are the $\sigma_1$- and $\sigma_2$-curvatures introduced by Viaclovsky \cite{Vi} and studied by many mathematicians.  

For metrics with $R>0$, Ge and Wang~\cite{GW} studied the following quotient
\begin{equation}\label{eq:intro-A}
A(M,g)\coloneqq\frac{\vol(g)\int_M\sigma_2\dv}{\bigl(\int_M\sigma_1\dv\bigr)^2}
\end{equation}
and proved $A(M,g)\le\frac13$, with equality exactly for metrics locally isometric to the round sphere. Later, Ge, Li and Lian~\cite{GLL} defined
\begin{equation}\label{eq:intro-A-sup}
A(M)\coloneqq\sup\bigl\{A(M,g):R_g>0\bigr\}
\end{equation}
and proved that $M$ is a spherical space form if $A(M)>-\frac{23}{384}$. The same conclusion holds if $A(M,g)>-\frac13$ for a metric of constant scalar curvature, or if $A(M,g)>-1$ for a Yamabe metric~\cite[Thms.~1.1--1.3 and Cor.~1.5]{GLL}. For the product metric on $\mathbb S^2\times \mathbb S^1$, the quotient \eqref{eq:intro-A} equals $-1$. They conjectured that $M$ is a spherical space form if $A(M)>-1$, and that $M$ is finitely covered by $\mathbb S^2\times \mathbb S^1$ if $A(M)=-1$~\cite[Conj.~1.6]{GLL}. If this conjecture were true, it would imply the inequality in the following question.

\smallskip \noindent{\bf Question 2.}
{\it If $M^3$ is not a spherical space form, does every metric with $R>0$ satisfy
\begin{equation}
   \vol(g)\int_M\sigma_2\dv\le -\bigl(\int_M\sigma_1\dv\bigr)^2,
\end{equation}
with equality if and only if $(M,g)$ is locally isometric to $\Sph^2\times \Sph^1$?
}

\smallskip

Unfornately, Questions~1 and~2 both have negative answers in dimension three (see Example~\ref{counterexampleforA}), and in dimensions $n\ge5$ even the weaker Question~3 below has a negative answer (see Example~\ref{prop:ngeq5}). We therefore turn to consider

\smallskip

\noindent{\bf Question 3.}
{\it Let $M^n$ be a closed manifold with $b_1(M)>0$. Does every metric with $R>0$ satisfy
\begin{equation}
  \vol(g)\int_M |\Ricz|^2\dv \ge \frac{1}{n(n-1)}\left( \int_M R\dv_g\right)^2,
\end{equation}
with equality only if $g$ is conformal to a metric locally isometric to $\Sph^{n-1}(a)\times \R$ for some $a>0$?
}

\medskip
Clearly, Question~3 is weaker than Questions~1 and~2. In particular, Gursky's theorem implies an affirmative answer to Question~3 in dimension $n=4$, and Theorem~\ref{thm:A} below gives an affirmative answer in dimension $n=3$, with equality only for $\mathbb S^2(a)\times\R$.

 For dimensions three and four, we actually  obtain a stronger result, under a weaker condition that $\int_M R_g dv_g>0$.

\begin{theorem}\label{thm:A}
Let $(M^3,g)$ be a closed connected  three-manifold with 
\begin{equation}\label{R>0}
\int_MR_g dv_g>0.
\end{equation}
I If $M$ is not diffeomorphic to a spherical space form $\mathbb S^3/\Gamma$,  then
\begin{equation}\label{eq:main}
\vol(g)\int_M|\Ricz|^2\dv
\ge \frac16\Bigl(\int_M R\dv\Bigr)^2.
\end{equation}
Equality holds if and only if $(M,g)$ is locally isometric to $\mathbb S^2(a)\times\R$ for $a>0$.
In particular, 
if 
\begin{equation*}
\vol(g)\int_M|\Ricz|^2\dv
<\frac16\Bigl(\int_M R\dv\Bigr)^2,
\end{equation*}
then $M$ is topologically a  spherical quotient.
\end{theorem}

In terms of $\sigma_1$ and $\sigma_2$ curvatures,  we have 
\begin{corollary}\label{coro_add} Let $(M^3,g)$ be a closed connected  three-manifold with \eqref{R>0}.
If 
\begin{equation}\label{eq:wish}
    \vol(g)\int_M\sigma_2\dv\ge - \Bigl(\int_M\sigma_1\dv\Bigr)^2 + \frac 13   \Bigl(\vol(g)  \int_M\sigma_1^2\dv-\Bigl(\int_M\sigma_1\dv\Bigr)^2 \Bigr),
\end{equation}
then $M^3$ is a spherical quotient or locally isometric to $\Sph^2 (a) \times \R$ for some $a>0$.
\end{corollary}
Equality \eqref{eq:wish} is the corrected inequality we took a long time to find.

\begin{theorem}\label{thm:total-scalar}
Let $(M^4,g)$ be a closed connected Riemannian four-manifold satisfying \eqref{R>0}.
 If some finite cover of $M$ has positive first Betti
number, then
\begin{equation}\label{eq:total-scalar}
\vol(g)\int_M| \Ricz_g|^2\dv_g
\ge \frac{1}{12}\left(\int_M R_g\dv_g\right)^2.
\end{equation}
Equality holds if and only if $R_g$ is constant
and $g$ is conformal to a metric locally isometric to $\mathbb S^3(a)\times\R$
for some $a>0$.
 
\end{theorem}

The above results are optimal in all aspects. It is clear that it is not true on sphercial forms and 
assumption \eqref{R>0} cannot be removed: a closed hyperbolic manifold with $b_1>0$ has $\Ricz=0$ and $\int_MR\dv<0$, so \eqref{eq:main} and \eqref{eq:total-scalar} fail for this example.

In the dimension $n=4$ case, if  
$g$ has positive Yamabe constant $Y(M,g)>0$, Gursky obtained the stronger inequality \eqref{eq:1.1}; other proofs were given by Bour and Carron~\cite{BC2} and are discussed in Remark~\ref{rem:compare-4d}. We give a new proof in Section 
\ref{sec:sec5}. When $g$ satisfies
only  \eqref{R>0}, Theorem \ref{thm:total-scalar} only shows 
\eqref{eq:total-scalar}. In this case the stronger inequality \eqref{eq:1.1} is false in general; see Remark~\ref{rem:L2-false}.

As consequences, Theorem \ref{thm:A} contains some part of the integral pinching theorem of Bour--Carron \cite[Thm.~5.6]{BC} and the results of Ge--Li--Lian for constant scalar curvature and Yamabe metrics \cite[Thms.~1.2, 1.3]{GLL}, with the constant $-\frac13$ of \cite[Thm.~1.2]{GLL} improved to $-1$. It does not imply \cite[Thm.~1.1]{GLL}, which concerns metrics with nonconstant scalar curvature.
\begin{corollary} \label{coro1}
   Let $(M,g)$ be a closed connected Riemannian three-manifold with $R_g>0$, and suppose $M$ is not a spherical space form. 
\begin{enumerate}[label=\textup{(\roman*)}] 
\item  If
\begin{equation}\label{eq:intro-BC}
\int_M|\Ricz|^2\dv\le \frac16\,Y(M,[g])^{3/2}(\min_M R)^{1/2},
\end{equation}
then $(M,g)$ is locally isometric to $\mathbb S^2(a)\times\R$.
\item If $g$ has constant scalar curvature, then $A(M,g)\le-1$.
Equality holds if and only if $(M,g)$ is locally isometric to $\mathbb S^2(a)\times\R$.
\end{enumerate}
\end{corollary}

\smallskip

\noindent{\it Comparison with Bour--Carron.}
Bour and Carron~\cite{BC} pass to a finite cover with $b_1>0$ and apply the Bochner formula to a harmonic one-form. Their argument implies that the Schr\"odinger operator $-\Delta+\frac12\rho_1$ cannot be positive, where $\rho_1(x)$ is the smallest eigenvalue of $\Ric$ at $x$, and combining this with the Yamabe--Sobolev inequality yields the pinching condition \eqref{eq:intro-BC}. After dividing \eqref{eq:main} by $\vol(g)$, the resulting lower bound matches \eqref{eq:intro-BC}  when $g$ is a Yamabe metric.

Our approach also begins with a finite cover with $b_1>0$. A key ingredient is the square identity \eqref{eq:n-square} established in Proposition~\ref{prop:n-square}. In dimension $n=3$, this identity is reduced to
\begin{align} 
\int|\Ricz|^2-\frac{(\int R_g dv_g)^2}{6\vol(g)}
=&\int\left|\Ricz-(\nabla\eta)^\circ
+\tfrac12(\eta\otimes\eta)^\circ\right|^2\\&+\frac16\int\left(\dvg\eta+|\eta|^2-\frac{\int R_g dv_g}{\vol(g)}\right)^2,\label{intro-square}
\end{align}
for a closed one-form $\eta$ solving the  following type equation
\begin{equation}\label{eq:riccati}
\tfrac52\dvg\eta+|\eta|^2=R.
\end{equation}
If such an $\eta$ exists,  \eqref{intro-square} implies the desired inequality \eqref{eq:main}.

Therefore, it remains to prove the existence of  solutions to equation \eqref{eq:riccati}. (For the rigidity statement we require two distinct solutions; see below.) Here the hypothesis $b_1>0$ plays a key role. Let $\alpha$ be a harmonic representative of a primitive
integral class in $H^1(M;\mathbb Z)$, and we seek a solution to \eqref{eq:riccati} of the form
\begin{equation}\label{eq:intro-ansatz}
\eta=\tfrac52\bigl(t\alpha+d\log v\bigr),
\qquad v>0.
\end{equation}
Substituting \eqref{eq:intro-ansatz} reduces \eqref{eq:riccati} to the linear equation
\begin{equation}\label{eq:3-Lt.0}
L_tv:=-\Delta v-2t\langle\alpha,dv\rangle
+\left(\frac{4}{25}R-t^2|\alpha|^2\right)v=0.
\end{equation}
A  solution to \eqref{eq:riccati} exists  when the principal eigenvalue $\lambda(L_t)$ of $L_t$ is zero. As $|t|\to\infty$ $\lambda(L_t)$ tends to $-\infty$, while under $R>0$, $\lambda(L_0)$ is positive. Hence there exists $t_0>0$ such that $\lambda(L_{t_0})=0$. This yields a solution of \eqref{eq:riccati} and thus proves \eqref{eq:main}, provided that $R_g>0$. A refinement of the same scheme leads to Theorem~\ref{thm:total-scalar}, which only assumes $\int_M R\dv>0$.

For the equality case, we start with   $\lambda (-t)=\lambda (t)$,  since $L_{-t}$ is the adjoint of $L_t$. Thus, we can obtain two forms $\eta_\pm$ of \eqref{eq:riccati} with $[\eta_+]=-[\eta_-]$.  Let $\xi=\frac12(\eta_+-\eta_-)$,  we show that $\xi^\sharp/|\xi|^2$ is Killing and that $M$ fibers over a circle. Passing to the associated cover, the lifted metric takes the form $h+b^2\,ds^2$, and we conclude that $b$ is constant and $h$ has constant positive Gauss curvature.

Under the weaker condition $\int_MR\dv>0$, the potential $\frac4{25}R$ may change sign, and $\lambda(0)$ may be negative. We then instead solve the shifted equation $$(\frac52+\tau)\dvg\eta+(1+\tau)|\eta|^2=R+\tau\bar R,$$ where $\bar R=\vol(g)^{-1}\int_MR\dv$ and $\tau\ge0$ is chosen so large such that $R+\tau\bar R>0$. The corresponding square identity (Lemma~\ref{lem:total-scalar-squares}) has the same structure as \eqref{intro-square}, with the coefficient $\frac16$ of the last square replaced by $\frac{1+2\tau}6$.

In the  $n=4$
case, the same strategy leads to a new proof of Gursky's result in $Y(M,g)>0$. The modified one described in the previous paragraph leads to  Theorem \ref{thm:total-scalar}.

We remark that the argument relies only on linear elliptic theory for the family $L_t$, in contrast with the fully nonlinear methods of \cite{CGY, GLW, CD, GW, GLL}. For a different approach to three-manifolds with $b_1>0$ and scalar curvature, through the level sets of harmonic maps to $\mathbb S^1$, see Stern~\cite{Stern2022}.

\medskip

\noindent{\it Organization of the rest of the paper.}
In Section~\ref{sec:prelim} we fix notations and carry out the finite-cover reduction. In Section~\ref{sec:dims} we construct the closed one-forms $\eta$ satisfying equation as \eqref{eq:riccati}, and prove the square identity. To illustrate the method, in Sections~4 and~\ref{sec:sec5} we prove Theorem~\ref{thm:A} under $R>0$ and give another proof of Gursky's theorem under $Y(M,[g])>0$. In Section~\ref{sec:total-scalar} we prove Theorems~\ref{thm:A} and~\ref{thm:total-scalar} under \eqref{R>0}.  The appendix contains the counterexamples we mentioned above.

\section{Preliminary}\label{sec:prelim}
Let $(M^n,g)$ be a closed, connected manifold of dimension $n\ge 3$.
For a symmetric $2$-tensor $S$, let
$S^\circ:=S-\frac1n(\tr S)g$ denote its trace-free part.
We identify one-forms and vector fields using the metric, and we use the conventions
\begin{equation}\label{eq:prelim-forms-notation}
(\nabla\eta)_{ij}=\nabla_i\eta_j,
\qquad
\dvg\eta=g^{ij}\nabla_i\eta_j,
\qquad
\langle\alpha,\beta\rangle=g^{ij}\alpha_i\beta_j.
\end{equation}

Unless otherwise indicated, trace-free parts are taken in dimension $n$.
In particular, if $\eta$ is a closed $1$-form, then
\[
(\nabla\eta)^\circ
=\nabla\eta-\frac{\dvg\eta}{n}g,
\qquad
(\eta\otimes\eta)^\circ
=\eta\otimes\eta-\frac{|\eta|^2}{n}g.
\]
Since $d\eta=0$, the tensor $\nabla\eta$ is symmetric, and
\begin{equation}\label{eq:n-square-decomposition}
|\nabla\eta|^2
=|(\nabla\eta)^\circ|^2+\frac1n(\dvg\eta)^2.
\end{equation}
We also set $\Ricz:=\Ric-\frac Rn g$.

By the contracted Bianchi identity, $\dvg\Ric=\frac12 dR$, it holds that
\begin{equation}\label{eq:n-divE}
\dvg \Ricz=\frac{n-2}{2n}\,dR.
\end{equation}

The quotient \eqref{eq:intro-A} and the quotient $\vol(g)\int_M|\Ricz|^2\dv/(\int_MR\dv)^2$ are invariant under rescaling and under finite Riemannian coverings.
If $\pi\colon(\hat M,\hat g)\to(M,g)$ is a finite Riemannian covering of degree $d$, then $\vol(\hat g)=d\,\vol(g)$ and each curvature integral appearing in these quotients is multiplied by $d$.

\begin{lemma}\label{lem:cover}
Let $M$ be a closed, connected $3$-manifold which is not a spherical space form.
Then $M$ admits a connected, orientable finite cover $\hat M$ with $b_1(\hat M)>0$.
In particular, $\hat M$ is not a spherical space form.
\end{lemma}

\begin{proof}
If $M$ is orientable,  Thurston's elliptization conjecture, which was proved by Perelman \cite{Pe,Pe2,Pe3}, states that a closed 3-manifold with finite fundamental group is spherical.
Since $M$ is not a spherical space form, it follows that $\pi_1(M)$ is infinite.
If $M$ is nonorientable, then $\chi(M)=0$ and $b_3(M)=0$, hence $b_1(M)=1+b_2(M)>0$ and $\pi_1(M)$ is again infinite.

The covering construction in \cite[Sect.~2, proof of Theorem~1, pp.~4--5]{GallowayLing2026},
using the virtual first Betti number result for aspherical summands \cite[Sect.~6, (I.2)]{AFW},
produces a connected, orientable finite cover $\hat M\to M$ with $H_2(\hat M;\Z)\neq 0$.
By Poincar\'e duality and the universal coefficient theorem,
\[
H_2(\hat M;\Z)\cong H^1(\hat M;\Z)\cong \Z^{b_1(\hat M)},
\]
hence $b_1(\hat M)>0$.
In particular, $\pi_1(\hat M)$ is infinite, which rules out spherical space forms.
\end{proof}

\begin{remark}\label{rem:psc}
If $M$ admits a metric of positive scalar curvature, then by the classification of $3$-manifolds with positive scalar curvature (see Schoen--Yau \cite{SY}, Gromov--Lawson \cite{GL,GL2}, Perelman \cite{Pe}), the orientation cover $M'$ is a connected sum of spherical space forms and copies of $\mathbb S^2\times \mathbb S^1$.
It is also standard that if $Y(M,[g])>0$, then $M'$ admits a metric of positive scalar curvature.
\end{remark}
\begin{lemma}\label{lem:descent}
Let $\pi\colon\hat M\to M$ be a finite Riemannian covering of degree $d$, and let $\hat g=\pi^*g$.
Then \eqref{eq:main} holds for $(M,g)$ if and only if it holds for $(\hat M,\hat g)$, and equality holds for one metric if and only if it holds for the other.
Moreover, for any $a>0$, $(M,g)$ is locally isometric to $\mathbb S^2(a)\times\R$ if and only if $(\hat M,\hat g)$ is.
\end{lemma}

\begin{proof}
Passing to the cover multiplies both sides of \eqref{eq:main} by $d^2$.
The local isometry statement follows since $\pi$ is a surjective local isometry.
\end{proof}

By Lemmas~\ref{lem:cover} and~\ref{lem:descent}, when proving Theorem~\ref{thm:A} we may assume that $M$ is orientable and $b_1(M)>0$.

The following lemma is a standard consequence of Tischler's fibration theorem
\cite{Tischler1970} and the standard properties
of Killing fields; see, for example,
\cite[Section~8.1]{Petersen}.

\begin{lemma}\label{lem:structure}
Let $(M,g)$ be a closed connected Riemannian manifold.
Let $\xi$ be a smooth closed nowhere-zero one-form
with period group $\ell\mathbb Z$, where $\ell>0$.
Suppose that
\(
X=\frac{\xi^\sharp}{|\xi|^2}
\)
is a Killing field for $g$.
Then the map
\[
\sigma(p)=\int_{p_0}^p\xi\pmod{\ell}
\]
defines a fiber bundle over $\mathbb R/\ell\mathbb Z$
with closed connected fiber $\Sigma$.
The connected infinite cyclic cover determined by
$\ker\sigma_*$, equipped with the lifted metric $\tilde g$, is isometric to
\[
(\Sigma\times\mathbb R, h+b(x)^2\,ds^2),
\]
where $h$ is a Riemannian metric on $\Sigma$ and $b$ is a smooth positive
function on $\Sigma$.

If $X$ has constant length, then $b$ is constant.
In particular, if $|X|_g=1$, the covering metric is
$h+ds^2$.
\end{lemma}

In Lemma~\ref{lem:structure}, $\partial_s$ is the lift of $X$.
The metric $h+b(x)^2\,ds^2$ on $\Sigma\times S^1_\ell$
and $(M,g)$ are locally isometric, since their common
covering space is $(\Sigma\times\R,h+b(x)^2\,ds^2)$.
Moreover, $\Sigma\times[0,\ell)$ is a common fundamental
domain for the deck action $\gamma(x, s)=(\psi(x), s+\ell)$, where $\psi$ is an isometry of $(\Sigma,h)$ satisfying
$b\circ\psi=b$.
Consequently, $(M,g)$ and
$(\Sigma\times S^1_\ell,h+b(x)^2\,ds^2)$  have the same volume and
the same integrals of scalar curvature and squared Ricci curvature.

\medskip

Now let $L$ be a scalar elliptic operator of the form
\begin{equation}\label{eq:forms-operator}
Lu=-\Delta_g u+2\langle\beta,du\rangle_g+q\,u,
\end{equation}
where $\beta$ is a smooth real one-form and $q$ is a smooth real-valued function.
Its adjoint with respect to $g$ is
\begin{equation}\label{eq:forms-adjoint}
L^*u=-\Delta_g u-2\langle\beta,du\rangle_g+(q-2\dvg\beta)\,u,
\end{equation}
where $\dvg\beta$ denotes the divergence of $\beta^\sharp$.

We recall the following standard facts (see \cite{KR}); (b) follows from $\int_M\bigl(\frac{L\varphi}{\varphi}-\lambda(L)\bigr)\varphi v^*\dv=\int_M\varphi\,(L^*v^*-\lambda(L)v^*)\dv=0$ and $\varphi v^*>0$.
\begin{lemma}\label{lem:principal}
There exists a real number $\lambda(L)$, called the principal eigenvalue, and smooth positive functions $v,v^*$ such that
\begin{equation}\label{eq:forms-eigenfunctions}
Lv=\lambda(L)v,\qquad L^*v^*=\lambda(L)v^*.
\end{equation}
\begin{enumerate}[label=\textup{(\alph*)}]
\item Every smooth positive eigenfunction of $L$ is a constant multiple of $v$, and the corresponding eigenvalue equals $\lambda(L)$.
\item For every smooth positive function $\varphi$ one has
\begin{equation}\label{eq:forms-eigenvalue-bound}
\min_M\frac{L\varphi}{\varphi}\le\lambda(L)\le\max_M\frac{L\varphi}{\varphi}.
\end{equation}
\end{enumerate}
\end{lemma}

\section{Special 1-forms and an integral identity}\label{sec:dims}

Throughout this section, $(M^n,g)$ is a closed, connected Riemannian manifold of dimension $n\ge 3$. We set
\[
\Ricz\coloneqq\Ric-\frac{R}{n}g,
\qquad
\bar R\coloneqq\frac1{\vol(g)}\int_M R\,dv_g.
\]
Unless otherwise indicated, all trace-free parts are taken in dimension $n$, and all integrals are over $M$ with respect to $dv_g$.

Assume for the moment that $M$ is orientable, that $b_1(M)>0$, and that
\begin{equation}\label{eq:n-spectral-hypothesis}
\lambda_1\left(-\Delta+\frac{2(n-1)}{(n+2)^2}R\right)>0.
\end{equation}
(In particular, this holds whenever $R>0$.)

Fix a primitive class in $H^1(M;\Z)$ and let $\alpha$ be its harmonic representative, so
\begin{equation}\label{eq:n-forms-harmonic}
d\alpha=0,\qquad \dvg\alpha=0,\qquad
\left\{\int_\gamma\alpha:\gamma\in\pi_1(M)\right\}=\Z.
\end{equation}
For $t\in\R$, consider the elliptic operator
\begin{equation}\label{eq:n-Lt}
L_tv=-\Delta v-2t\langle\alpha,dv\rangle
+\left(\frac{2(n-1)}{(n+2)^2}R-t^2|\alpha|^2\right)v.
\end{equation}
Write $\lambda(t)\coloneqq\lambda(L_t)$ and let $v_t>0$ be the (normalized) positive principal eigenfunction, chosen so that $\max_M v_t=1$.

\begin{lemma}\label{lem:n-lambda}
Assume $\dvg\alpha=0$ and $\alpha \neq 0$. Under the assumption  \eqref{eq:n-spectral-hypothesis}, we have 
\begin{enumerate}[label=\textup{(\alph*)}]
\item $\lambda\colon\R\to\R$ is continuous and even.
\item For every $t\in\R$,
\begin{equation}\label{eq:n-forms-lambda-integral}
\lambda(t)\vol(g)
=\frac{2(n-1)}{(n+2)^2}\int R
-t^2\int|\alpha|^2-\int|d\log v_t|^2.
\end{equation}
\item There exists $t_0>0$ such that $\lambda(t_0)=0$.
\end{enumerate}
\end{lemma}

\begin{proof}
Since $\dvg\alpha=0$, the adjoint satisfies $L_t^*=L_{-t}$, hence $\lambda(-t)=\lambda(t)$.

For continuity, fix $s\in\R$ and use $v_s>0$ as a test function in \eqref{eq:forms-eigenvalue-bound}. A direct computation gives
\begin{equation}\label{eq:n-forms-continuity}
\frac{L_t v_s}{v_s}
=\lambda(s)-2(t-s)\langle\alpha,d\log v_s\rangle-(t^2-s^2)|\alpha|^2
\longrightarrow\lambda(s)
\end{equation}
uniformly on $M$ as $t\to s$, which implies $\lambda(t)\to\lambda(s)$, by \eqref{eq:forms-eigenvalue-bound}.

To prove \eqref{eq:n-forms-lambda-integral}, by the definition of $\lambda(t)$, we have
\begin{align*}
\lambda(t)\vol(g)
&=\int\frac{L_tv_t}{v_t}\\
&=\frac{2(n-1)}{(n+2)^2}\int R-t^2\int|\alpha|^2-\int|d\log v_t|^2,
\end{align*}
using $\Delta v_t/v_t=\Delta\log v_t+|d\log v_t|^2$, together with $\int\Delta\log v_t=0$ and $\int\langle\alpha,d\log v_t\rangle=0$.

By \eqref{eq:n-spectral-hypothesis} we have $\lambda(0)>0$. Since $\alpha\not\equiv 0$, \eqref{eq:n-forms-lambda-integral} shows that $\lambda(t)\to-\infty$ as $|t|\to\infty$. By continuity, there exists some $t_0>0$ such that $\lambda(t_0)=0$.
\end{proof}

\begin{proposition}\label{prop:n-forms}
Assume that $M$ is orientable and $b_1(M)>0$, and let $\alpha$ be as above. Let $t_0>0$ be as in Lemma~\ref{lem:n-lambda}.
Then the one-forms
\begin{equation}\label{eq:n-forms-eta-pair}
\eta_\pm\coloneqq\frac{n+2}{n-1}\left(\pm t_0\alpha+d\log v_{\pm t_0}\right)
\end{equation}
are smooth and closed, and represent the nontrivial classes $[\eta_\pm]=\pm\frac{n+2}{n-1}t_0[\alpha]\neq 0$. Moreover, each $\eta_\pm$ satisfies
\begin{equation}\label{eq:n-eta}
\frac{n+2}{2}\dvg\eta+\frac{n-1}{2}|\eta|^2=R.
\end{equation}
\end{proposition}

\begin{proof}
From \eqref{eq:n-forms-harmonic} and \eqref{eq:n-forms-eta-pair} we immediately obtain $d\eta_\pm=0$ and $[\eta_\pm]=\pm\frac{n+2}{n-1}t_0[\alpha]\neq 0$.

Next,
\begin{equation*}
\frac{n+2}{2}\dvg\eta_\pm
=\frac{(n+2)^2}{2(n-1)}\left(\frac{\Delta v_{\pm t_0}}{v_{\pm t_0}}-|d\log v_{\pm t_0}|^2\right)
\end{equation*}
and
\begin{equation*}
\frac{n-1}{2}|\eta_\pm|^2
=\frac{(n+2)^2}{2(n-1)}\left(t_0^2|\alpha|^2+|d\log v_{\pm t_0}|^2\pm 2t_0\langle\alpha,d\log v_{\pm t_0}\rangle\right).
\end{equation*}
Adding these identities and using $L_{\pm t_0}v_{\pm t_0}=0$ yields \eqref{eq:n-eta}.
\end{proof}

We now record the basic identities involving $\Ricz$ and a closed one-form $\eta$ that will be used to factor the square appearing in Proposition~\ref{prop:n-square}.

\begin{lemma}\label{lem:ndim}
For every smooth closed one-form $\eta$ on $M$ one has
\begin{align}
\int\langle \Ricz,(\nabla\eta)^\circ\rangle
&=\frac{n-2}{2n}\int R\,\dvg\eta,
\label{eq:n-EH}\\
\int\langle \Ricz,(\eta\otimes\eta)^\circ\rangle
&=\int\left(\frac{n-1}{n}(\dvg\eta)^2-|(\nabla\eta)^\circ|^2-\frac1nR|\eta|^2\right),
\label{eq:n-EB}\\
|(\eta\otimes\eta)^\circ|^2
&=\frac{n-1}{n}|\eta|^4,
\label{eq:n-HB0}\\
\int\langle(\nabla\eta)^\circ,(\eta\otimes\eta)^\circ\rangle
&=-\frac{n+2}{2n}\int(\dvg\eta)|\eta|^2.
\label{eq:n-HB}
\end{align}
\end{lemma}

\begin{proof} These identities are well known; we include the short proofs.
Since $\Ricz$ is trace-free, integration by parts and \eqref{eq:n-divE} yields
\begin{align*}
\int\langle \Ricz,(\nabla\eta)^\circ\rangle
&=\int\langle \Ricz,\nabla\eta\rangle
=-\int\langle\dvg \Ricz,\eta\rangle \\
&=-\frac{n-2}{2n}\int\langle dR,\eta\rangle
=\frac{n-2}{2n}\int R\,\dvg\eta,
\end{align*}
which proves \eqref{eq:n-EH}.

For \eqref{eq:n-EB}, integrating the Bochner--Weitzenb\"ock formula for $\eta$ and using $d\eta=0$ and \eqref{eq:n-square-decomposition}, we obtain
\begin{align}
\int\Ric(\eta,\eta)
&=\int\left((\dvg\eta)^2-|\nabla\eta|^2\right)\notag\\
&=\int\left(\frac{n-1}{n}(\dvg\eta)^2-|(\nabla\eta)^\circ|^2\right).
\label{eq:n-square-bochner}
\end{align}
Since $\langle \Ricz,(\eta\otimes\eta)^\circ\rangle=\Ric(\eta,\eta)-\frac1nR|\eta|^2$, this yields \eqref{eq:n-EB}.

Identity \eqref{eq:n-HB0} follows from $|(\eta\otimes\eta)^\circ|^2=|\eta\otimes\eta|^2-\frac1n|\eta|^4$.
Finally,
\begin{align*}
\langle(\nabla\eta)^\circ,(\eta\otimes\eta)^\circ\rangle
&=\langle\nabla\eta,\eta\otimes\eta\rangle-\frac1n(\dvg\eta)|\eta|^2
=\frac12\langle\eta,d|\eta|^2\rangle-\frac1n(\dvg\eta)|\eta|^2,
\end{align*}
which implies  \eqref{eq:n-HB}, by integration by parts.
\end{proof}

We now prove a key identity.

\begin{proposition}\label{prop:n-square}

For every smooth closed one-form $\eta$ one has
\begin{align}
&\int\left|\Ricz-(\nabla\eta)^\circ+\frac12(\eta\otimes\eta)^\circ\right|^2\notag\\
=&\int|\Ricz|^2-\frac{n-1}{4n}\int\bigl((n-2)\dvg\eta+|\eta|^2\bigr)^2
+\frac{(n-3)(n-4)}4\int (\dvg\eta)^2\notag\\
&\qquad+\int \frac1n\bigl((n-2)\dvg\eta+|\eta|^2\bigr)
\left(\frac{n+2}{2}\dvg\eta+\frac{n-1}{2}|\eta|^2-R\right).
\label{eq:expandn}
\end{align}
In particular, if $\eta$ satisfies \eqref{eq:n-eta}, then
\begin{equation}\label{eq:n-square-mean}
\int\bigl((n-2)\dvg\eta+|\eta|^2\bigr)=\frac{2\vol(g)}{n-1}\bar R,
\end{equation}
and
\begin{align}
\int|\Ricz|^2-\frac{\vol(g)\bar R^2}{n(n-1)}
&=\int\left|\Ricz-(\nabla\eta)^\circ+\frac12(\eta\otimes\eta)^\circ\right|^2-\frac{(n-3)(n-4)}4\int(\dvg\eta)^2
\notag\\
&\quad+\frac{n-1}{4n}\int\left((n-2)\dvg\eta+|\eta|^2-\frac{2\bar R}{n-1}\right)^2
\label{eq:n-square}
\end{align}
\end{proposition}

\begin{proof}
Expanding the square, the term $\int|(\nabla\eta)^\circ|^2$ cancels against the corresponding term of \eqref{eq:n-EB}, and \eqref{eq:n-EH}--\eqref{eq:n-HB} show that   
\begin{align}
&\int\left|\Ricz-(\nabla\eta)^\circ+\frac12(\eta\otimes\eta)^\circ\right|^2\notag\\
=&\int|\Ricz|^2+\int\bigl(\frac{n-1}n(\dvg\eta)^2+\frac{n+2}{2n}\dvg\eta\,|\eta|^2+\frac{n-1}{4n}|\eta|^4-\frac{n-2}nR\dvg\eta-\frac1nR|\eta|^2\bigr).
\end{align}
Then, collecting the coefficients  and
rearranging the above formula of the five terms, we have
\begin{align}
&\int\left|\Ricz-(\nabla\eta)^\circ+\frac12(\eta\otimes\eta)^\circ\right|^2\notag\\
=&\int|\Ricz|^2-\frac{n-1}{4n}\int\bigl((n-2)\dvg\eta+|\eta|^2\bigr)^2
+\frac{(n-3)(n-4)}4\int (\dvg\eta)^2\notag\\
&\qquad+\int \frac1n\bigl((n-2)\dvg\eta+|\eta|^2\bigr)
\left(\frac{n+2}{2}\dvg\eta+\frac{n-1}{2}|\eta|^2-R\right),
\label{eq:n-square-factorization}
\end{align}
which is \eqref{eq:expandn}.

Assume now that \eqref{eq:n-eta} holds. Since $\int\dvg\eta=0$, we compute
\begin{equation*}
\int\bigl((n-2)\dvg\eta+|\eta|^2\bigr)
=\frac{2}{n-1}\int\left(\frac{n+2}{2}\dvg\eta+\frac{n-1}{2}|\eta|^2\right)
=\frac{2\vol(g)}{n-1}\bar R,
\end{equation*}
which is \eqref{eq:n-square-mean}. Moreover,
\begin{align*}
\int\bigl((n-2)\dvg\eta+|\eta|^2\bigr)^2
=\int\left((n-2)\dvg\eta+|\eta|^2-\frac{2\bar R}{n-1}\right)^2+\frac{4\vol(g)\bar R^2}{(n-1)^2}.
\end{align*}
Substituting into \eqref{eq:expandn} yields \eqref{eq:n-square}.
\end{proof}

For $n=3,4$ the tensor $\Ricz-(\nabla\eta)^\circ+\frac12(\eta\otimes\eta)^\circ$ in \eqref{eq:n-square} has a geometric meaning. For $n=3$ it is the traceless Ricci tensor of the torsion-free connection $D_XY=\nabla_XY+\frac12\eta(X)Y+\frac12\eta(Y)X$, which is projectively equivalent to the Levi-Civita connection; if $\eta=-df$ locally, the Ricci tensor of $D$ is the $1$-Bakry--\'Emery Ricci tensor $\Ric+\nabla^2f+\frac12df\otimes df$, studied by Wylie and Yeroshkin~\cite{WylieYeroshkin}. For $n=4$ it is the traceless Ricci tensor of the closed Weyl connection $D_XY=\nabla_XY+\frac12\bigl(\eta(X)Y+\eta(Y)X-g(X,Y)\eta^\sharp\bigr)$, which is locally the Levi-Civita connection of $e^{2f}g$ when $\eta=2df$.

For $n\ge5$, the term $-\frac{(n-3)(n-4)}4\int(\dvg\eta)^2$ in \eqref{eq:n-square} is nonpositive, so the method used in the next sections (for $n=3,4$) cannot directly extend. In Example~\ref{prop:ngeq5} below we construct metrics in these dimensions for which the corresponding pinching bound fails.

\section{3-dimensional case: Proof of Theorem \ref{thm:A} under the stronger assumption $R_g>0$}
Throughout, we work on an orientable closed $3$-manifold with $b_1(M)>0$, and  we prove Theorem~\ref{thm:A} under the stronger assumption $R_g>0$.

\medskip

Taking $\eta=\eta_\pm$ in \eqref{eq:n-forms-eta-pair}, the identity \eqref{eq:n-square} becomes
\begin{align}
\int|\Ricz|^2-\frac{\vol(g)\bar R^2}{6}
=&\int\left|\Ricz-(\nabla\eta)^\circ
+\tfrac12(\eta\otimes\eta)^\circ\right|^2+\frac16\int\left(\dvg\eta+|\eta|^2-\bar R\right)^2.
\label{eq:square}
\end{align}
Since $R>0$ implies \eqref{eq:n-spectral-hypothesis}, Proposition~\ref{prop:n-forms} yields a solution of \eqref{eq:riccati}; inserting it into \eqref{eq:square} yields \eqref{eq:main}.

It remains to analyze the equality case. Assume equality holds in \eqref{eq:main}. For each of the two forms $\eta\in\{\eta_+,\eta_-\}$ from Proposition~\ref{prop:n-forms} (with $n=3$), both squares in \eqref{eq:square} vanish, i.e.
\begin{equation}\label{oneformrelationwith R}
\dvg\eta+|\eta|^2-\bar R=0
\end{equation}
 and
\begin{equation}\label{einstein-equality}
\Ricz=(\nabla\eta)^\circ-\frac12(\eta\otimes\eta)^\circ.
\end{equation}
Combining \eqref{oneformrelationwith R} with \eqref{eq:riccati}, we obtain
\begin{equation}\label{eq:eqcase}
R=\tfrac52\bar R-\tfrac32|\eta|^2.
\end{equation}
Equivalently,
\begin{align}\label{eq:einstein-another-form}
\Ric&=\nabla\eta-\tfrac12\eta\otimes\eta+
\Bigl(\tfrac R3-\tfrac{\dvg\eta}3+\tfrac{|\eta|^2}6\Bigr)g\notag\\
&=\nabla\eta-\tfrac12\eta\otimes\eta+\tfrac12\bar R\,g,
\end{align}
where the last equality uses \eqref{oneformrelationwith R} and \eqref{eq:eqcase}. In particular, \eqref{eq:eqcase} implies $|\eta_+|=|\eta_-|$ pointwise.

Define
\begin{equation}\label{eq:equality-sum-difference}
\xi\coloneqq\tfrac12(\eta_+-\eta_-),
\qquad
\vartheta\coloneqq\tfrac12(\eta_++\eta_-).
\end{equation}
It is easy to see
\begin{equation}\label{eq:equality-xi-class}
\langle\xi,\vartheta\rangle
=\tfrac14(|\eta_+|^2-|\eta_-|^2)=0,
\qquad
[\xi]=\tfrac52t_0[\alpha]\ne0.
\end{equation}

\begin{lemma}\label{lem:xi}
The forms $\xi$ and $\vartheta$ satisfy
\begin{equation}\label{eq:equality-xi-derivatives}
\nabla\xi=\tfrac12(\vartheta\otimes\xi+\xi\otimes\vartheta),
\qquad
d|\xi|^2=|\xi|^2\vartheta.
\end{equation}
Moreover, $\xi\neq 0$ on $M$, and $X\coloneqq\xi^\sharp/|\xi|^2$ is a Killing field with $\xi(X)=1$.
\end{lemma}

\begin{proof}
Subtracting \eqref{eq:einstein-another-form} for $\eta_+$ and $\eta_-$, we have
\[
2\nabla\xi=\frac12\bigl(\eta_{+} \otimes \eta_{+}-\eta_{-} \otimes \eta_{-}\bigr).
\]
Then writing $\eta_\pm=\vartheta\pm\xi$ yields
\begin{equation}\label{equality-1}
2\nabla\xi
=\tfrac12\bigl((\vartheta+\xi)\otimes(\vartheta+\xi)-(\vartheta-\xi)\otimes(\vartheta-\xi)\bigr)
=\vartheta\otimes\xi+\xi\otimes\vartheta.
\end{equation}
Contracting \eqref{equality-1} with $\xi$ and using $\langle\vartheta,\xi\rangle=0$, it holds that
\[
\xi^j\nabla_i\xi_j=\frac12|\xi|^2\vartheta_i.
\]
Thus,  it holds that $d|\xi|^2=|\xi|^2\vartheta$. Along any smooth path $\gamma$,
\[
|\xi|^2(\gamma(1))=|\xi|^2(\gamma(0))\exp\Bigl(\int_\gamma\vartheta\Bigr).
\]
By connectedness, $|\xi|^2$ is either identically zero or positive everywhere; since $[\xi]\neq 0$, we have $|\xi|^2>0$.

Now $X^\flat=\xi/|\xi|^2$, so differentiating and using \eqref{equality-1}, we get
\begin{align}
\nabla_i X_j
&=\frac{\nabla_i\xi_j}{|\xi|^2}-\frac{\xi_j\,\partial_i |\xi|^2}{|\xi|^4}
 =\frac{\xi_i\vartheta_j-\vartheta_i\xi_j}{2|\xi|^2}.
\label{eq:equality-killing}
\end{align}
The skew-symmetry of \eqref{eq:equality-killing} shows that $X$ is Killing. By construction, $\xi(X)=1$.
\end{proof}

Due to Lemma \ref{lem:structure}, we now show that any metric of the form $h+b(x)^2\,ds^2$ on $\Sigma\times \mathbb S^1_\ell$ that attains equality in \eqref{eq:main} must be a round product.

\begin{proposition}\label{prop:warped}
Let $N=\Sigma\times \mathbb S^1_\ell$ carry the warped product metric $h+b(x)^2\,ds^2$, where $(\Sigma,h)$ is a closed connected surface with Gauss curvature $K$ and $b$ is a smooth positive function on $\Sigma$.
Assume equality holds in \eqref{eq:main} and $\int_N R\,dv_g>0$. Then there exists a positive constant $b$ such that
$$(N, g) \cong\left\{\begin{array}{l}\left(\mathbb{S}^2 \times S_{\ell}^1, {2}/{\bar{R}} g_{\mathbb{S}^2}+b^2 d s^2\right) \\ \left(\mathbb{R P}^2 \times S_{\ell}^1, {2}/{\bar{R}} g_{\mathbb{R P}^2}+b^2 d s^2\right).\end{array}\right.$$

\end{proposition}

\begin{proof}
Let $\nu=b^{-1}\partial_s$. The warped product curvature formulas yield
\begin{align}\label{eq:warped1}
\Ric(\nu,\nu) &=-\frac{\Delta_h b}{b},\\
\label{eq:warped2}
    \Ric(\nu,\cdot)|_{T\Sigma}&=0,\\
\label{eq:warped3}
    \Ric|_{T\Sigma}&=Kh-\frac{\nabla_h^2 b}{b},\\
\label{eq:warped4}
     R&=2K-\frac{2\Delta_h b}{b}.
\end{align}
Since $dv=b\,dA_h\,ds$,  where $dA_h$ is the area form w.r.t. $h$, \eqref{eq:warped1} implies
\[
\int_N\Ric(\nu,\nu)\dv=-\ell\int_\Sigma\Delta_h b\,dA_h=0.
\]
Set $P=\nu^\flat\otimes\nu^\flat-\frac13g$. It is easy to check that $\langle \Ricz,P\rangle=\Ric(\nu,\nu)-\frac R3$, and
$|P|^2_g=\frac23$. Combining with equality in
\eqref{eq:main}, we obtain
\begin{align}
\int_N\Bigl|\Ricz+\frac{\bar R}{2}P\Bigr|^2\dv
&=\int_N|\Ricz|^2\dv+\bar R\int_N\Ric(\nu,\nu)\dv
  -\frac{\bar R}{3}\int_N R\dv+\frac{vol(N)\bar R^2}{6}\nonumber\\
&=\int_N|\Ricz|^2\dv-\frac{vol(N)\bar R^2}{6}=0, 
\label{eq:equality-warped-square}
\end{align}
which yields \begin{equation}\label{rigidity1}
    \Ricz=-\frac{\bar R}{2}P
\end{equation}
 It is clear that $P|_{T\Sigma}=-\frac 13 h$, and hence $\left.\Ricz\right|_{T \Sigma}=\frac{\bar{R}}{6} h$.  By \eqref{eq:warped3}, we know that $\left.\Ricz\right|_{T \Sigma}=\left(K-\frac{R}{3}\right) h-\frac{\nabla_h^2 b}{b}$. Thus, 
\begin{equation}
    \nabla_h^2 b=b\left(K-\frac{R}{3}-\frac{\bar R}{6}\right) h, 
\end{equation}
that is, 
\begin{equation}\label{eq:twoeq1}
\nabla_h^2 b=\tfrac12(\Delta_h b)\,h.
\end{equation}
Hence $\nabla_h b$ is a conformal vector field on $(\Sigma, h)$.
Moreover, $\Ricz(\nu,\nu)=-\frac13\bar R$ by \eqref{rigidity1}.
Together with \eqref{eq:warped1} and \eqref{eq:warped4}, we have
\begin{equation}\label{eq:twoeq}
2Kb+\Delta_h b=\bar R\,b.
\end{equation}
By  \eqref{eq:twoeq1}, we know that  on $\Sigma$
\begin{align}
\frac12 d(\Delta_h b)
&=\operatorname{div}_h\!\left(\frac12(\Delta_h b)h\right)
 =\operatorname{div}_h(\nabla_h^2 b)\notag\\
&=d(\Delta_h b)+\operatorname{Ric}_h(\nabla_h b,\cdot)\notag=d(\Delta_h b)+K\,db,\notag
\end{align}
and thus $d\Delta_h b=-2K\,db$. 

Differentiating \eqref{eq:twoeq} then yields $2b\,dK=\bar R\,db$.
Since $\nabla_h b$ is a conformal Killing field, by the Kazdan--Warner
identity \cite{KW}, which holds for conformal vector fields on every closed surface \cite{BE}, we have
\[
0=\int_\Sigma\langle\nabla_h K,\nabla_h b\rangle_h\,dA_h
=\frac{\bar R}{2}\int_\Sigma\frac{|\nabla_h b|_h^2}{b}\,dA_h.
\]
Thus $b$ is constant, and \eqref{eq:twoeq} implies $K=\bar R/2$.
\end{proof}

\begin{remark}
The structure  for the proof of the equality case is similar to the one given  in Lafontaine's theorem~\cite{Lafontaine2009} that a compact Riemannian three-manifold carrying two independent solutions of the static equation $\nabla^2f-(\Delta f)g-f\Ric=0$ is isometric to a standard $\mathbb S^1\times\mathbb S^2$, $\mathbb S^1\times\mathbb{RP}^2$ or $\mathbb S^3$: a Killing field is built from two solutions, the metric is reduced to a warped product $h+u^2\,ds^2$, and the warping function satisfies $\nabla_h^2u=\frac12(\Delta_hu)h$ on the surface.
\end{remark}

\begin{proof}[Proof of the equality case in Theorem~\ref{thm:A} under $R>0$]
If $(M,g)$ is locally isometric to $\mathbb S^2(a)\times\R$, the Ricci
eigenvalues are $a^{-2},a^{-2},0$. Thus $R=2a^{-2}$ is constant and
$|\Ricz|^2=\frac23a^{-4}=\frac16R^2$, which yields equality in
\eqref{eq:main}.

Conversely, assume equality in \eqref{eq:main}. By
Lemmas~\ref{lem:cover} and~\ref{lem:descent}, we may assume $M$ is
orientable and $b_1(M)>0$. By Lemma~\ref{lem:structure}, the metric $g$
is locally isometric to $h+b(x)^2\,ds^2$ on
$\Sigma\times \mathbb S^1_\ell$ and attains equality in \eqref{eq:main}.

The total
scalar curvature of $h+b(x)^2\,ds^2$ is positive because the total
scalar curvature of $(M,g)$ is positive. Proposition~\ref{prop:warped}
implies that $b$ is constant and the Gauss curvature of  surface metric $h$
is  constant  $\frac12\bar R>0$. The original metric
is therefore locally isometric to $\mathbb S^2(a)\times\R$ with
$a^{-2}=\frac12\bar R$.
\end{proof}

From the previous proof, one can see that 
\begin{equation}\label{eq:consequences-spectral-hypothesis}
\lambda_1\bigl(-\Delta+\tfrac4{25}R\bigr)>0
\end{equation}
is enough. In fact we will improve it to $\int_M R_g dv_g>0$.

\begin{proof}
[\noindent{\it Proof of Corollary \ref{coro1}.}]
For (i), we use $u\equiv1$ as a test function in the definition of the Yamabe
constant. Thus $Y(M,[g])\le vol(g)^{-1/3}\int R$. Combining the bound with
$\min_M R\le\bar R$, we obtain
\begin{equation}\label{eq:consequences-yamabe-bound}
Y(M,[g])^{3/2}(\min_M R)^{1/2}\le \frac{(\int R)^2}{vol(g)}.
\end{equation}
Hypothesis \eqref{eq:intro-BC} and \eqref{eq:consequences-yamabe-bound} imply
$vol(g)\int|\Ricz|^2\le\frac16(\int R)^2$. 

Together with Theorem~\ref{thm:A}, $\vol(g)\int|\Ricz|^2=\frac16(\int R)^2$, and the equality case of Theorem~\ref{thm:A} gives the conclusion.

For (ii), since $\sigma_1=\frac R4$ and $\sigma_2=\frac{R^2}{48}-\frac12|\Ricz|^2$ in dimension three, $$A(M,g)=\frac{\vol(g)\int R^2}{3(\int R)^2}-\frac{8\vol(g)\int|\Ricz|^2}{(\int R)^2}.$$ If $R$ is constant, the first term equals $\frac13$ and, by Theorem~\ref{thm:A}, the second is at least $\frac43$, with equality exactly when equality holds in \eqref{eq:main}. 
\end{proof}

(i) is a result of Bour and Carron~\cite[Thm.~5.7 and Rem.~5.8]{BC}.
(ii) includes the Yamabe metric case of \cite[Thm.~1.3]{GLL} and sharpens the constant $-\frac13$ of \cite[Thm.~1.2]{GLL} to $-1$.

\section{Four-dimensional case: A new proof of Gursky's Theorem }
\label{sec:sec5}

The method from the previous section also applies in dimension $n=4$. Let $Y(M,[g])$ denote the Yamabe constant. In dimension four, the condition $Y(M,[g])>0$ is equivalent to
$\lambda_1\bigl(-\Delta+\frac16R\bigr)>0$.

We restate the pinching theorem mentioned in the Introduction.

\begin{theorem}[Gursky~\cite{Gu}]\label{thm:four-dimensional}
Let $(M^4,g)$ be closed and connected with $Y(M,[g])>0$. Suppose that $M$ admits a connected finite cover with $b_1>0$. Then
\begin{equation}\label{eq:strong4}
\int_M|\Ricz|^2\,dv_g\ge\frac1{12}\int_M R^2\,dv_g,
\end{equation}
and hence
\begin{equation}\label{eq:main4}
\vol(g)\int_M|\Ricz|^2\,dv_g
\ge\frac1{12}\left(\int_MR\,dv_g\right)^2.
\end{equation}
Equality in \eqref{eq:strong4} holds if and only if $g$ is conformal to a metric locally isometric to $\mathbb S^3(a)\times\R$ for some $a>0$. Equality in \eqref{eq:main4} holds if and only if, in addition, $R_g$ is constant.
\end{theorem}

The inequality \eqref{eq:strong4} was proved by Gursky~\cite{Gu}, and later reproved by Bour--Carron~\cite{BC2}. We provide another proof, following the same idea as in the three-dimensional case. In dimension four, $\frac{2(n-1)}{(n+2)^2}=\frac16$, so $L_0$ is the conformal Laplacian, and the form $\eta$ constructed below is the Lee form of a closed Weyl structure with vanishing conformal scalar curvature. From this point of view, the square identity \eqref{eq:square4} is the integrated identity for closed Weyl structures of Calderbank--Pedersen~\cite[Sect.~9]{CP} (see also \cite{PPS}), and the new point here is the construction of $\eta$ through the family $L_t$.

\begin{proposition}\label{prop:4dforms}
Let $(M^4,g)$ be closed, connected, and orientable, with $Y(M,[g])>0$ and $b_1(M)>0$. Let $\alpha$ be the harmonic representative of a primitive integral class, i.e. $\alpha$ satisfies \eqref{eq:n-forms-harmonic}.
Then there exist $t_0>0$ and positive functions $v_{\pm t_0}$ such that
\begin{equation}\label{eq:4-forms-eta-pair}
\eta_\pm\coloneqq 2\left(\pm t_0\alpha+d\log v_{\pm t_0}\right)
\end{equation}
are closed, satisfy $[\eta_\pm]=\pm2t_0[\alpha]\ne0$, and solve
\begin{equation}\label{eq:riccati4}
3\dvg\eta+\frac32|\eta|^2=R.
\end{equation}
\end{proposition}

\begin{proof}
This is Proposition~\ref{prop:n-forms} with $n=4$, since
$\frac{2(n-1)}{(n+2)^2}=\frac16$ and
$\frac{n+2}{n-1}=2$.
\end{proof}

For every smooth closed one-form $\eta$ in dimension four, \eqref{eq:n-square-factorization} is
\begin{align}
\int\left|\Ricz-(\nabla\eta)^\circ+\frac12(\eta\otimes\eta)^\circ\right|^2\notag
=&\int|\Ricz|^2-\frac3{16}\int(2\dvg\eta+|\eta|^2)^2\notag\\
&+\frac14\int(2\dvg\eta+|\eta|^2)
\left(3\dvg\eta+\frac32|\eta|^2-R\right).
\label{eq:Q4}
\end{align}
If $\eta$ satisfies \eqref{eq:riccati4}, that is, $2\dvg\eta+|\eta|^2=\frac23R$, then 
\begin{equation}\label{eq:square4}
\int|\Ricz|^2-\frac1{12}\int R^2
=\int\left|\Ricz-(\nabla\eta)^\circ+\frac12(\eta\otimes\eta)^\circ\right|^2.
\end{equation}
Equivalently,
\begin{align}
\int|\Ricz|^2-\frac{\vol(g)\bar R^2}{12}
&=\int\left|\Ricz-(\nabla\eta)^\circ+\frac12(\eta\otimes\eta)^\circ\right|^2+\frac1{12}\int(R-\bar R)^2.
\label{eq:4-square-mean}
\end{align}

\begin{proof}[Proof of the inequalities in Theorem~\ref{thm:four-dimensional}]
Passing to a connected finite cover, we may assume $b_1(M)>0$; if necessary we further pass to the orientation cover. Positivity of the conformal Laplacian (equivalently $Y(M,[g])>0$) is preserved under finite covers: the pullback of a positive principal eigenfunction remains a positive eigenfunction with the same eigenvalue, and Lemma~\ref{lem:principal}(a) shows it is still principal.

Proposition~\ref{prop:4dforms} then provides a closed one-form satisfying \eqref{eq:riccati4}. The nonnegativity of the right-hand sides of \eqref{eq:square4} and \eqref{eq:4-square-mean} yields \eqref{eq:strong4} and \eqref{eq:main4} on the cover.

If $\pi\colon(\widetilde M,\widetilde g)\to(M,g)$ is a finite Riemannian cover of degree $d$, then both sides of \eqref{eq:strong4} scale by $d$, while both sides of \eqref{eq:main4} scale by $d^2$. Hence the inequalities descend to $M$.
\end{proof}

\begin{remark}
In dimension four, \eqref{eq:square4} gives the stronger inequality \eqref{eq:strong4}. By contrast, in dimension three we only have the averaged inequality \eqref{eq:main} coming from \eqref{eq:square}; an analog of \eqref{eq:strong4} for $n=3$ is false (see the final section). The reason is: in $n=4$, the term $\int (  2\dvg\eta+|\eta|^2)^2$ in \eqref{eq:Q4} gives 
$\int R^2$, while in $n=3$ the corresponding term gives $\vol(g)\bar R^2=\vol(g)^{-1}(\int R)^2$.
\end{remark}

\begin{remark}\label{rem:4dconformal}
In dimension four, the Schouten tensor is $A_g=\frac12\bigl(\Ric-\frac R6g\bigr)$, and
\begin{equation}\label{eq:4-conformal-deficit}
\int_M\left(|\Ricz|^2-\frac1{12}R^2\right)\,dv_g
=-8\int_M\sigma_2\bigl(g^{-1}A_g\bigr)\,dv_g.
\end{equation}
Hence \eqref{eq:strong4} is equivalent to
\[
\int_M \sigma_2\bigl(g^{-1}A_g\bigr)\,dv_g\le 0,
\]
a conformally invariant condition.

By the Chern--Gauss--Bonnet formula,
\begin{equation}\label{GBC}
8\pi^2\chi(M)=\frac14\int_M|W_g|^2\,dv_g+4\int_M\sigma_2\bigl(g^{-1}A_g\bigr)\,dv_g,
\end{equation}
so \eqref{eq:strong4} yields the estimate $\int_M|W_g|^2\,dv_g\ge 32\pi^2\chi(M)$. If \(M\) is nonorientable, we obtain the same estimates by applying the estimate to its orientation double cover equipped with the pullback metric.

The forms $\eta_\pm$ of Proposition~\ref{prop:4dforms} define closed Weyl structures with vanishing conformal scalar curvature: \eqref{eq:riccati4} is exactly this condition. The set of such structures was studied by Higa~\cite{Higa2001}. Assuming the set is nonempty, he showed that each structure is determined by its cohomology class, and that the set of classes is symmetric under $u\mapsto-u$. Proposition~\ref{prop:4dforms} shows that the set is nonempty whenever $Y(M,[g])>0$ and $b_1(M)>0$.
\end{remark}

\subsection{Equality case for four-dimensional case}We now turn to the equality case. Assume equality holds in \eqref{eq:strong4}. We first treat the case where $M$ is orientable and $b_1(M)>0$; the general case will follow by a finite-cover argument.

Let $\eta\in\{\eta_+,\eta_-\}$ be the two closed one-forms given by Proposition~\ref{prop:4dforms}. Since the right-hand side of \eqref{eq:square4} is a nonnegative integral, equality in \eqref{eq:strong4} forces the square to vanish, i.e.
\begin{equation}\label{eq:4-einstein-equality}
\Ricz=(\nabla\eta)^\circ-\frac12(\eta\otimes\eta)^\circ.
\end{equation}
Together with \eqref{eq:riccati4}, this is equivalent to
\begin{align}
\Ric
&=\nabla\eta-\frac12\eta\otimes\eta+\left(\frac R4-\frac{\dvg\eta}{4}+\frac{|\eta|^2}{8}\right)g\notag\\
&=\nabla\eta-\frac12\eta\otimes\eta+\frac12\bigl(\dvg\eta+|\eta|^2\bigr)g.
\label{eq:4-ricci-equality}
\end{align}
Unlike in dimension three, these identities do not force $|\eta_+|=|\eta_-|$.

As before, set
\begin{equation}\label{eq:4-sum-difference}
\xi\coloneqq\frac12(\eta_+-\eta_-),
\qquad
\vartheta\coloneqq\frac12(\eta_++\eta_-).
\end{equation}
Then $\xi$ and $\vartheta$ are closed and $[\xi]=2t_0[\alpha]\ne0$. Subtracting \eqref{eq:riccati4} for $\eta_+$ and $\eta_-$ yields
\begin{equation}\label{eq:4-xi-divergence}
\dvg\xi=-\langle\xi,\vartheta\rangle,
\end{equation}
which is not  necessarily zero, different to the case $n=3$.

\begin{lemma}\label{lem:4-xi}
The forms $\xi$ and $\vartheta$ satisfy
\begin{equation}\label{eq:4-xi-derivatives}
\nabla\xi=\frac12\left(
\vartheta\otimes\xi+\xi\otimes\vartheta
-\langle\xi,\vartheta\rangle g\right),
\qquad
d|\xi|^2=|\xi|^2\vartheta.
\end{equation}
The form $\xi\neq 0$ for every $x\in M$. Moreover,
\begin{equation}\label{eq:4-conformal-metric}
\widehat g=|\xi|^2g,
\qquad
X=\frac{\xi^\sharp}{|\xi|^2}
\end{equation}
define a smooth metric and a unit parallel vector field for
$\widehat g$, with $\widehat g(X,\cdot)=\xi$ and $\xi(X)=1$.
Here $\xi^\sharp$ and $|\xi|$ in \eqref{eq:4-conformal-metric} are
computed with respect to $g$.
\end{lemma}

\begin{proof}
Subtracting \eqref{eq:4-ricci-equality} for $\eta_+$ and $\eta_-$ gives
\begin{align*}
2\nabla\xi
&=\vartheta\otimes\xi+\xi\otimes\vartheta
-\bigl(\dvg\xi+2\langle\xi,\vartheta\rangle\bigr)g\\
&=\vartheta\otimes\xi+\xi\otimes\vartheta
-\langle\xi,\vartheta\rangle g,
\end{align*}
where we used \eqref{eq:4-xi-divergence}. Contracting with $\xi^j$,
we obtain
\begin{equation*}
\xi^j\nabla_i\xi_j
=\frac12\left(\vartheta_i|\xi|^2
+\xi_i\langle\xi,\vartheta\rangle
-\xi_i\langle\xi,\vartheta\rangle\right)
=\frac12|\xi|^2\vartheta_i.
\end{equation*}
 Since $\xi^j\nabla_i\xi_j=\frac12|\xi|^2\vartheta_i$, we have $d|\xi|^2=|\xi|^2\,\vartheta$. Hence $|\xi|^2$ is either identically zero or strictly positive on each connected component. Because $[\xi]\ne0$, we conclude $|\xi|^2>0$ everywhere.

Since $X^\flat=\xi/|\xi|^2$ with respect to $g$, we compute
\begin{align}
\nabla_iX_j
=\frac{\nabla_i\xi_j}{|\xi|^2}
-\frac{\xi_j\partial_i |\xi|^2}{|\xi|^4}
=\frac{\xi_i\vartheta_j-\vartheta_i\xi_j}{2|\xi|^2}
-\frac{\langle\xi,\vartheta\rangle}{2|\xi|^2}g_{ij}.
\label{eq:4-conformal-killing}
\end{align}
Therefore,
\begin{align}\label{eq:b}
\mathcal L_X g=-\frac{\langle\xi,\vartheta\rangle}{|\xi|^2}\,g.
\end{align}
Thus $X$ is a conformal Killing field for $g$. Now consider the conformal metric $\widehat g=|\xi|^2\,g$ (which agrees with \eqref{eq:4-conformal-metric}).

Using \eqref{eq:4-xi-derivatives} we have $X(|\xi|^2)=\langle\xi,\vartheta\rangle$, and together with \eqref{eq:b} this implies
\begin{align}
\mathcal L_X\widehat g
&=X(|\xi|^2)g+|\xi|^2\mathcal L_Xg=0.
\end{align}
Moreover, $\widehat g(X,X)=1$ and $\widehat g(X,\cdot)=\xi$.
Since $X$ is Killing for $\widehat g$, the tensor
$\widehat\nabla\xi$ is skew-symmetric. Since $d\xi=0$, it is also
symmetric. Hence $\widehat\nabla\xi=0$, so $X$ is parallel.
\end{proof}

By Lemma \ref{lem:structure}, we know $\hat g=h+\ ds^2$, where
$(\Sigma,h)$ is a closed connected three-manifold.
Now let us give the splitting theorem for the four-dimensional manifolds. 
\begin{proposition}\label{prop:4-product}
Let $N=\Sigma\times \mathbb S^1_\ell$ with $\widehat g=h+ds^2$, where
$(\Sigma,h)$ is a closed connected three-manifold.
Suppose equality holds in \eqref{eq:strong4} and
$\int_N R_{\widehat g}\,dv_{\widehat g}>0$.
Then $h$ has constant positive sectional curvature.
In particular, $(N,\widehat g)$ is locally isometric to
$\mathbb S^3(a)\times\R$ for some $a>0$.
\end{proposition}

\begin{proof}
In this proof all curvature tensors refer to $\widehat g$. Let $\nu\coloneqq\partial_s$. For  $\widehat g=h+ds^2$,  we have
\begin{equation}\label{eq:4-product-curvature}
\Ric(\nu,\nu)=0,
\qquad
\Ric|_{T\Sigma}=\Ric_h,
\qquad
R=R_h.
\end{equation}
Set $P\coloneqq\nu^\flat\otimes\nu^\flat-\frac14\widehat g$. Then $\langle\Ricz,P\rangle=-\frac14R$ and $|P|^2=\frac34$. 

Since the equality holds in \eqref{eq:strong4},  we obtain
\begin{align}
\int_N\left|\Ricz+\frac R3P\right|^2\,dv_{\widehat g}
&=\int_N\left(|\Ricz|^2+\frac{2R}{3}\langle\Ricz,P\rangle+\frac{R^2}{9}|P|^2\right)\,dv_{\widehat g}\notag\\
&=\int_N|\Ricz|^2\,dv_{\widehat g}-\frac1{12}\int_N R^2\,dv_{\widehat g}=0,
\label{eq:4-product-square}
\end{align}
by the conformal invariant as in Remark \ref{rem:4dconformal},
and then $\Ricz=-\frac R3P$ pointwise.

Restricting to $T\Sigma$ (where $P=-\frac14h$) yields
\[
\Ricz|_{T\Sigma}=\Ric_h-\frac{R_h}{4}h=\frac{R_h}{12}h,
\]
that is, $\Ric_h=\frac{R_h}{3}h$, an Einstein metric.   Since $\int_N R\,dv_{\widehat g}>0$ and $R=R_h$, we have $R_h>0$,  hence $h$ has constant positive sectional curvature. (In dimension three, an Einstein metric has constant sectional curvature and  $R_h$ is constant by Schur's lemma).
\end{proof}

\begin{proof}[Proof of the equality statements in Theorem~\ref{thm:four-dimensional}]
The ``if'' direction is easy.

For the converse, assume first that $M$ is orientable and $b_1(M)>0$. By Lemma~\ref{lem:structure} and~\ref{prop:4-product}, $\widehat g=|\xi|^2g$ is locally isometric to $\mathbb S^3(a)\times\R$.

It remains to show that $\widehat g$ descends to a metric
in the conformal class $[g]$ on $M$.  Choose a connected finite normal cover $\pi\colon\widehat M\to M$ which dominates both the orientation cover and a cover with $b_1(\widehat M)>0$. And we have  a metric $\widehat g=e^{2u}\pi^*g$ on $\widehat M$ locally isometric to $\mathbb S^3(a)\times\R$.

Every deck transformation of $\pi$ is conformal for $\widehat g$. Lifting to the universal cover, we obtain a conformal self-map of the round cylinder $\mathbb S^3(a)\times\R$, which can be identified with
\[
\bigl(\R^4\setminus\{0\},\ a^2|x|^{-2}|dx|^2\bigr).
\]
By Liouville's theorem, such a conformal self-map is the restriction of a M\"obius transformation preserving $\{0,\infty\}$, hence it can only be  the form $x\mapsto\lambda Ax$ or $x\mapsto\lambda Ax/|x|^2$, where $\lambda>0$ and $A\in O(4)$. In both cases the metric $a^2|x|^{-2}|dx|^2$ is preserved, so each deck transformation is in fact an isometry of $\widehat g$. Consequently, $u$ is deck-invariant and $\widehat g$ descends to a metric in the conformal class $[g]$ on $M$.

Finally,
\begin{align*}
\vol(g)\int|\Ricz|^2-\frac1{12}\left(\int R\right)^2
&=\vol(g)\left(\int|\Ricz|^2-\frac1{12}\int R^2\right)+\frac{\vol(g)}{12}\int(R-\bar R)^2.
\end{align*}
Both terms are nonnegative. Hence equality in \eqref{eq:main4} is equivalent to equality in \eqref{eq:strong4} together with constant scalar curvature, which proves both equality statements.
\end{proof}
\begin{remark}
\label{rem:compare-4d}Gursky~\cite{Gu} and Bour--Carron~\cite{BC2} proved \eqref{eq:strong4} under the hypothesis $b_1(M)>0$. Sung~\cite[Theorem~1.3(iii)]{Su} obtained the same inequality for closed oriented four-manifolds with $Y(M,[g])\ge0$ and finite covers of arbitrarily large degree. For $Y(M,[g])>0$ and infinite $\pi_1(M)$, \eqref{eq:strong4} also follows from Chang--Gursky--Yang~\cite[Theorem~A and Corollary~B]{CGY} together with the Bonnet--Myers theorem.
\end{remark}

In dimension three, equality in \eqref{eq:main} forces $g$ to be locally isometric to $\mathbb S^2(a)\times\R$. In dimension four, however, equality in \eqref{eq:main4} does \emph{not} preclude a negative Ricci eigenvalue, as the following example shows.

\smallskip
\noindent\textit{Example.} Let $u$ be a positive, nonconstant periodic solution of $u''=u-u^3$ with period $L$, obtained by taking a small oscillation about $u=1$. The metric
\[
 g\coloneqq u(t)^2\,(dt^2+g_{\mathbb S^3})
\]
on $\mathbb S^3\times \mathbb S^1_L$ has constant scalar curvature $R_g=6u^{-3}(u-u'')=6$. Since $g$ is conformal to the round product, it attains equality in \eqref{eq:main4}. At a minimum of $u$, the Ricci eigenvalue in the circle direction equals $-3u''/u^3<0$, so $(\mathbb S^3\times \mathbb S^1_L,g)$ is not locally isometric to $\mathbb S^3(a)\times\R$.

\section{Pinching under positive total scalar curvature}\label{sec:total-scalar}

In this section we prove Theorems~\ref{thm:A} and~\ref{thm:total-scalar}, which require only the assumption $\int R>0$. 
For a smooth closed one-form $\eta$, set
\[
S_\eta=(n-2)\dvg\eta+|\eta|^2 .
\]
We first prove an integral identity.
\begin{lemma}\label{lem:total-scalar-squares}
Let $(M^n,g)$ be a closed Riemannian manifold of dimension $n\ge3$.
Suppose that $\eta$ satisfies
\begin{equation}\label{eq:total-scalar-shift}
\left(\frac{n+2}{2}+(n-2)\tau\right)\dvg\eta
+\left(\frac{n-1}{2}+\tau\right)|\eta|^2
=R+\frac{2\tau}{n-1}\bar R
\end{equation}
for some $\tau\ge0$. Then
\begin{align}
\int|\Ricz|^2-\frac{\vol(g)\bar R^2}{n(n-1)}
=&\int\left|\Ricz-(\nabla\eta)^\circ
+\tfrac12(\eta\otimes\eta)^\circ\right|^2+\frac{n-1+4\tau}{4n}
       \int\left(S_\eta-\frac{2\bar R}{n-1}\right)^2\notag\\
&\quad-\frac{(n-3)(n-4)}4\int(\dvg\eta)^2 .
\label{eq:total-scalar-squares-n}
\end{align}
In particular, if $n=3$ or $n=4$, then
\begin{align}
\int|\Ricz|^2-\frac{\vol(g)\bar R^2}{n(n-1)}
=\int\left|\Ricz-(\nabla\eta)^\circ
                 +\tfrac12(\eta\otimes\eta)^\circ\right|^2
+\frac{n-1+4\tau}{4n}
       \int\left(S_\eta-\frac{2\bar R}{n-1}\right)^2 .
\label{eq:total-scalar-squares}
\end{align}
\end{lemma}

\begin{proof}
As above, all integrals are over $M$ with respect to $dv_g$. Since $M$ is
closed, we have $\int\dvg\eta=0$. Integrating
\eqref{eq:total-scalar-shift}, we obtain
\begin{align*}
\int|\eta|^2
=\frac{2\vol(g)\bar R}{n-1},
\end{align*}
and 
\begin{align*}
\int S_\eta
=(n-2)\int\dvg\eta+\int|\eta|^2
=\frac{2\vol(g)\bar R}{n-1}.
\end{align*}
It follows that
\begin{align}\label{eq:c1}
\int\left(S_\eta-\frac{2\bar R}{n-1}\right)=0.
\end{align}

We next compute the  terms in 
\eqref{eq:expandn}. Rearranging \eqref{eq:total-scalar-shift} gives
\begin{align*}
R=\frac{n+2}{2}\dvg\eta+\frac{n-1}{2}|\eta|^2
 +\tau\left(S_\eta-\frac{2\bar R}{n-1}\right).
\end{align*}
Multiplying this expression by $S_\eta$ and using its definition, we
find
\begin{align*}
RS_\eta
&=\frac{(n+2)(n-2)}{2}(\dvg\eta)^2
  +\frac{n+2+(n-1)(n-2)}{2}(\dvg\eta)|\eta|^2\\
&\quad+\frac{n-1}{2}|\eta|^4
  +\tau S_\eta\left(S_\eta-\frac{2\bar R}{n-1}\right),
\end{align*}
and thus,
\begin{align*}
\text{LHS of}\,\, \eqref{eq:expandn} =&\int |\Ricz|^2+\int \bigg[ \left(\frac{n-1}{n}
             -\frac{(n+2)(n-2)}{2n}\right)(\dvg\eta)^2-\frac{(n-1)(n-2)}{2n}(\dvg\eta)|\eta|^2\\
       &\quad  -\frac{n-1}{4n}|\eta|^4-\frac{\tau}{n}S_\eta
               \left(S_\eta-\frac{2\bar R}{n-1}\right)\bigg].
\end{align*}
We also have
\begin{align*}
-\frac{n-1}{4n}S_\eta^2
&=-\frac{(n-1)(n-2)^2}{4n}(\dvg\eta)^2\\
&\quad-\frac{(n-1)(n-2)}{2n}(\dvg\eta)|\eta|^2
      -\frac{n-1}{4n}|\eta|^4.
\end{align*}
Combining these calculations with \eqref{eq:expandn}, we get
\begin{align}\label{middle term}
\int\left|\Ricz-(\nabla\eta)^\circ
             +\tfrac12(\eta\otimes\eta)^\circ\right|^2
&=\int|\Ricz|^2-\frac{n-1}{4n}\int S_\eta^2
  +\frac{(n-3)(n-4)}4\int(\dvg\eta)^2\\
&\quad-\frac{\tau}{n}\int S_\eta
              \left(S_\eta-\frac{2\bar R}{n-1}\right).
\end{align}
By \eqref{eq:c1}, it is easy to check 
\begin{align*}
\int S_\eta\left(S_\eta-\frac{2\bar R}{n-1}\right)
&=\int\left(S_\eta-\frac{2\bar R}{n-1}\right)^2
 +\frac{2\bar R}{n-1}\int\left(S_\eta-\frac{2\bar R}{n-1}\right)\\
&=\int\left(S_\eta-\frac{2\bar R}{n-1}\right)^2,
\end{align*}
and
\begin{align*}
\int S_\eta^2
&=\int\left(S_\eta-\frac{2\bar R}{n-1}
                        +\frac{2\bar R}{n-1}\right)^2\\
&=\int\left(S_\eta-\frac{2\bar R}{n-1}\right)^2
 +\frac{4\bar R}{n-1}\int\left(S_\eta-\frac{2\bar R}{n-1}\right)
 +\frac{4\vol(g)\bar R^2}{(n-1)^2}\\
&=\int\left(S_\eta-\frac{2\bar R}{n-1}\right)^2
 +\frac{4\vol(g)\bar R^2}{(n-1)^2}.
\end{align*}
Substituting these two formulas into \eqref{middle term}
yields
\begin{align*}
\int\left|\Ricz-(\nabla\eta)^\circ
+\tfrac12(\eta\otimes\eta)^\circ\right|^2
=&\int|\Ricz|^2+\frac{(n-3)(n-4)}4\int(\dvg\eta)^2\\
&-\frac{n-1+4\tau}{4n}
       \int\left(S_\eta-\frac{2\bar R}{n-1}\right)^2-\frac{\vol(g)\bar R^2}{n(n-1)}
,
\end{align*}
and we obtain \eqref{eq:total-scalar-squares}.
\end{proof}

\medskip
Now we begin to finish the whole proof.
\begin{proof}[Proof of Theorems~\ref{thm:A} and~\ref{thm:total-scalar}]
By Lemma~\ref{lem:cover} if $n=3$, and by assumption if $n=4$, $M$ has a connected finite cover with
$b_1>0$. We may also take its orientation cover, since pullback on
first real cohomology is injective for finite covers. Both sides of
\eqref{eq:main} and of \eqref{eq:total-scalar} are multiplied by the square of the degree,
and $\bar R$ is preserved.

\medskip\noindent
\textbf{Claim.}
\emph{As Proposition \ref{prop:n-forms}, Equation~\eqref{eq:total-scalar-shift} admits two  solutions
$\eta_\pm$ whose cohomology classes satisfy $[\eta_+]=-[\eta_-]$.}

\medskip

Let $\alpha$ be the harmonic representative of a primitive integral
class satisfying \eqref{eq:n-forms-harmonic}. 
Let
\[
P_\tau=\frac{R+\frac{2\tau}{n-1}\bar R}
             {(\frac{n-1}{2}+\tau)c_\tau^2},
\]
where \(c_\tau=\frac{\frac{n+2}{2}+(n-2)\tau}{\frac{n-1}{2}+\tau}\) and take $\tau>\max \left\{0,-\frac{n-1}{2 \bar{R}} \min _M R\right\}$.
As in Proposition~\ref{prop:n-forms}, we reduce
\eqref{eq:total-scalar-shift} to $L_tv=0$, where
\[
L_t v
=
-\Delta v-2t\langle\alpha,dv\rangle
+
\left(
\frac{R+\frac{2\tau}{n-1}\bar R}
     {\left(\frac{n-1}{2}+\tau\right)c_\tau^2}
-t^2|\alpha|^2
\right)v.
\]
Since $P_\tau>0$ for $\tau>\max \left\{0,-\frac{n-1}{2 \bar{R}} \min _M R\right\}$,
the argument of Lemma~\ref{lem:n-lambda} yields $t_0>0$ and positive
functions $v_{\pm t_0}$ satisfying $L_{\pm t_0}v_{\pm t_0}=0$.
The corresponding forms
\begin{equation}\label{two forms in generalcase}
\eta_\pm=c_\tau\bigl(\pm t_0\alpha+d\log v_{\pm t_0}\bigr)
\end{equation}
therefore solve \eqref{eq:total-scalar-shift}. Since $\alpha$ is closed
and $d\log v_{\pm t_0}$ is exact, the forms $\eta_\pm$ are closed and
\[
[\eta_\pm]=\pm c_\tau t_0[\alpha]\ne 0.
\]
This proves the claim.

\smallskip
Now applying \eqref{eq:total-scalar-squares} to
$\eta_\pm$, we obtain \eqref{eq:main} and \eqref{eq:total-scalar} respectively.

For the equality case, by \eqref{eq:total-scalar-squares} we have
\[
S_{\eta_\pm}=(n-2)\dvg\eta_\pm+|\eta_\pm|^2=\frac{2\bar R}{n-1},
\]
and
\[
\Ricz=(\nabla\eta_\pm)^\circ-\tfrac12(\eta_\pm\otimes\eta_\pm)^\circ .
\]

Together with  equation \eqref{eq:total-scalar-shift},  $\eta_\pm$ satisfy \eqref{eq:n-eta}, that is, \begin{equation}\label{eq:n-eta-1}
\frac{n+2}{2}\dvg\eta_\pm+\frac{n-1}{2}|\eta_\pm|^2=R.
\end{equation} 

For $n=3$ and $n=4$, we now arrive at the same equalities as in the case
$R>0$, and $\xi=\tfrac12(\eta_+-\eta_-)$ has period group
$c_\tau t_0\Z$. The same rigidity argument therefore applies.

\end{proof}

\begin{remark}
\label{rem:L2-false}
(i) The hypothesis $\int_MR_g\dv_g>0$ can hold on a manifold which admits no
metric of positive scalar curvature.  For $\tilde g=u^{4/(n-2)}g$,  one has $\int_MR_{\tilde g}\,dv_{\tilde g}=\int_M\bigl(\frac{4(n-1)}{n-2}|du|^2+Ru^2\bigr)\dv$, which can be  positive for suitable $u>0$ in the following way. Take $u=\varepsilon+\varphi_\delta$ with a smooth cut-off function $\varphi_\delta$ of radius $\delta$ such that  $\int|d\varphi_\delta|^2\sim\delta^{n-2}$ and $\int\varphi_\delta^2\sim\delta^n$.  Then $\int_MR_{\tilde g}\,dv_{\tilde g}>0$  for  small  $\delta$ and then choosing sufficiently small $\varepsilon$ such that $u$ is positive. For instance, Theorem~\ref{thm:A} applies to such metrics in the conformal class of a hyperbolic metric.

(ii) For every metric on the manifolds $M^4$ with $\chi(M)\le0$, by the Chern--Gauss--Bonnet formula \eqref{GBC}, we know that
 \eqref{eq:strong4} and \eqref{eq:total-scalar} hold. If $Y(M,[g])>0$, \eqref{eq:total-scalar} follows from Theorem~\ref{thm:four-dimensional}.  If $Y(M,[g])=0$, then $[g]$ contains a scalar-flat metric, for which \eqref{eq:strong4} is trivial; since the left-hand side of \eqref{eq:4-conformal-deficit} is conformally invariant, \eqref{eq:strong4} holds for every metric in $[g]$. 
 
 The new case of Theorem~\ref{thm:total-scalar} is therefore $\chi(M)>0$ and $Y(M,[g])<0$, in which \eqref{eq:strong4} can fail in the following way.
 By Millson~\cite{Mil}, there exists a closed hyperbolic
four-manifold $(M,h)$ with $b_1(M)>0$.
By part~(i) of Remark \ref{rem:L2-false}, the conformal class $[h]$ contains metric $g$
with positive total scalar curvature and  $W_g=0$,  and hence
\[
\int_M|\Ricz_g|^2\,dv_g-\frac1{12}\int_MR_g^2\,dv_g
=-16\pi^2\chi(M)<0.
\]

\end{remark}

\section{Appendix: Counterexamples}
In this section we will provide various examples related to Question 1, 2 and 3 as mentioned in the introduction.
\subsection{Counterexamples in dimensions at least five}

The negative term in \eqref{eq:n-square} is not merely an artifact
of the proof: the inequality of Question~3 (and hence that of Question~1) is false
for $n\ge5$.

\begin{example}\label{prop:ngeq5}
On $\mathbb S^{n-1}\times \mathbb S^1_L$, let
\[
g_0=dt^2+g_{\mathbb S^{n-1}},
\qquad
g_\varepsilon=e^{2\varepsilon\cos(\omega t)}g_0,
\qquad
\omega=\frac{2\pi}{L},
\]
where $g_{\mathbb S^{n-1}}$ has sectional curvature one.
Then
\begin{equation}\label{eq:n-counterexample-expansion}
\frac{n(n-1)\vol(g_\varepsilon)
\int|\Ricz_{g_\varepsilon}|^2\,dv_{g_\varepsilon}}
{\left(\int R_{g_\varepsilon}\,dv_{g_\varepsilon}\right)^2}
=1+\frac{\varepsilon^2}{2}
\left(\omega^4-2(n-2)\omega^2+4\right)+O(\varepsilon^3).
\end{equation}
For every $n\ge5$, choosing $\omega^2=n-2$ and sufficiently small
$|\varepsilon|\ne0$, we have $R_{g_\varepsilon}>0$ and
\begin{equation}\label{eq:n-counterexample}
\vol(g_\varepsilon)\int|\Ricz_{g_\varepsilon}|^2\,dv_{g_\varepsilon}
<\frac1{n(n-1)}
\left(\int R_{g_\varepsilon}\,dv_{g_\varepsilon}\right)^2.
\end{equation}
\end{example}

\begin{proof}
Set $u(t)=\varepsilon\cos(\omega t)$.
By computations, we know that 
\begin{align}
R_{g_\varepsilon}
&=(n-1)e^{-2u}
\left((n-2)-2u''-(n-2)(u')^2\right),\notag\\
|\Ricz_{g_\epsilon}|^2
&=\frac{(n-1)(n-2)^2}{n}e^{-4u}
\left(1+u''-(u')^2\right)^2,\notag\\
dv_{g_\varepsilon}&=e^{nu}\,dt\,dv_{\mathbb S^{n-1}}.
\label{eq:n-counterexample-curvatures}
\end{align}
Integrating the scalar curvature formula by parts, we obtain 
\[
\int R_{g_\varepsilon}\,dv_{g_\varepsilon}
=(n-1)(n-2)\vol(\mathbb S^{n-1})
\int_0^L e^{(n-2)u}\left(1+(u')^2\right)\,dt.
\]
Using the mean values of $\cos(\omega t)$, $\cos^2(\omega t)$,
and $\sin^2(\omega t)$, we obtain
\begin{align}
\frac{\vol(g_\varepsilon)}{\vol(g_0)}
&=1+\frac{n^2}{4}\varepsilon^2+O(\varepsilon^3),\notag\\
\frac{\int R_{g_\varepsilon}\,dv_{g_\varepsilon}}
{(n-1)(n-2)\vol(g_0)}
&=1+\left(\frac{(n-2)^2}{4}+\frac{\omega^2}{2}\right)
\varepsilon^2+O(\varepsilon^3),\notag\\
\frac{n\int|\Ricz_{g_\varepsilon}|^2\,dv_{g_\varepsilon}}
{(n-1)(n-2)^2\vol(g_0)}
&=1+\left(\frac{(n-4)^2}{4}-(n-3)\omega^2
+\frac{\omega^4}{2}\right)\varepsilon^2
+O(\varepsilon^3).
\label{eq:n-counterexample-integrals}
\end{align}
Thus, we have
\eqref{eq:n-counterexample-expansion}.

For $\omega^2=n-2$, the coefficient of $\varepsilon^2$ is
$\frac12(4-(n-2)^2)<0$ when $n\ge5$.
Thus \eqref{eq:n-counterexample} holds for sufficiently small
nonzero $\varepsilon$. Since
$R_{g_0}=(n-1)(n-2)>0$, scalar curvature remains positive for all
sufficiently small $|\varepsilon|$.
\end{proof}

For $n=3$ and $n=4$, the quadratic coefficient in
\eqref{eq:n-counterexample-expansion} is
$\frac12((\omega^2-1)^2+3)>0$ and $\frac12(\omega^2-2)^2\ge0$, respectively.

\subsection{Counterexamples to Questions 1 and 2 in dimension three}\label{sec:counter}

Let $(\Sigma,h)$ be a closed surface of Gauss curvature $K$.
For a  function $b>0$ on $\Sigma$, we equip
$\Sigma\times \mathbb S^1_L$ with the metric
\begin{equation}\label{eq:counter-warped-metric}
 g=h+b^2dt^2.
\end{equation}
Let  $\Delta$ and $dA_h$ denote the surface Laplacian and area form,
respectively. The trace-free Hessian on $\Sigma$ is denoted by $\hesz b$.
Recall that
\begin{equation}\label{eq:counter-curvature-relations}
 \sigma_2=\tfrac3{16}R^2-\tfrac12|\Ric|^2,\qquad
 |\Ricz|^2=|\Ric|^2-\tfrac13R^2.
\end{equation}

\begin{proposition}\label{prop:formulas}
For the metric in \eqref{eq:counter-warped-metric}, we have
\begin{align}
 \sigma_1 b
 &=\tfrac12(Kb-\Delta b), \label{eq:counter-sigma1}\\
 \sigma_2 b
 &=-\tfrac14K^2b-\tfrac12K\Delta b-\frac{|\hesz b|^2}{2b}, \label{eq:counter-sigma2}\\
 2|\Ricz|^2b
 &=\frac{(2Kb+\Delta b)^2}{3b}+\frac{2|\hesz b|^2}{b}, \label{eq:counter-ricci}
\end{align}
and
\begin{equation}\label{eq:warped-A-plus-one}
\begin{aligned}
&\frac{\left(\int_\Sigma Kb\,dA_h\right)^2}{\int_\Sigma b\,dA_h}\left(A(M,g)+1\right)\\
=&
\biggl(
2\int_\Sigma\langle\nabla_h K,\nabla_h b\rangle_h\,dA_h
-\int_\Sigma b(K-\bar K_b)^2\,dA_h
-2\int_\Sigma\frac{|(\nabla_h^2 b)^\circ|_h^2}{b}\,dA_h
\biggr),
\end{aligned}
\end{equation}
where $\bar{K}_b=\frac{\int_{\Sigma} K b d A_h}{\int_{\Sigma} b d A_h}$.

\end{proposition}

\begin{proof}
The equality \eqref{eq:warped4} yield \eqref{eq:counter-sigma1}.
From \eqref{eq:warped1}--\eqref{eq:warped3} we obtain
\begin{equation}\label{eq:counter-norms}
 \begin{aligned}
 |\Ric|^2&=(\Delta b/b)^2+|Kh-\nabla^2 b/b|^2,\\
 |\nabla^2 b|^2&=\tfrac12(\Delta b)^2+|\hesz b|^2.
 \end{aligned}
\end{equation}
Substituting \eqref{eq:counter-norms} and the formula for $R$ in
\eqref{eq:warped4} into \eqref{eq:counter-curvature-relations}, we obtain
\eqref{eq:counter-sigma2} and \eqref{eq:counter-ricci}. Finally,
we use \eqref{eq:counter-sigma1}, \eqref{eq:counter-sigma2}, and
$dv=b\,dA_h\,dt$ in \eqref{eq:intro-A}. 
Since $\int_\Sigma\Delta b\,dA_h=0$,
we obtain 
\begin{equation}\label{eq:Aformula}
\bigl(\int_\Sigma Kb\,dA_h\bigr)^2 A(M, g)=-
 \bigl(\int_\Sigma b\,dA\bigr)
\int_\Sigma\bigl(K^2b+2K\Delta b+2\frac{|\hesz b|^2}{b}\bigr)\,dA_h.
\end{equation}
Then integrating by parts, we have \eqref{eq:warped-A-plus-one}.
\end{proof}

For a rotationally symmetric metric on $\mathbb S^2$, we use the moment coordinate
$x\in[-1,1]$ and the angle $\theta\in\R/2\pi\Z$.
Let $F$ be smooth on $[-1,1]$, and set
\begin{equation}\label{eq:moment}
 \begin{aligned}
 h&=\frac{dx^2}{F(x)}+F(x)\,d\theta^2,
\end{aligned}
\end{equation}
where \begin{equation}\label{eq:c2}
 F(\pm1)=0,\qquad F'(\pm1)=\mp2,\qquad F>0\quad\text{on }(-1,1).\end{equation}
The conditions in \eqref{eq:c2} ensure smoothness at the poles. The area form is $dA=dx\,d\theta$, so the area is $4\pi$.
For a rotationally symmetric function $b=b(x)$, we compute
\begin{equation}\label{eq:counter-moment-curvature}
 K=-\tfrac12 F'',\qquad
 \Delta b=(Fb')',\qquad
 |\hesz b|^2=\tfrac12 F^2(b'')^2.
\end{equation}
The round metric corresponds to $F=1-x^2$.

\begin{example}\label{counterexampleforA}
On $\mathbb S^2\times \mathbb S^1$, we use the coordinates
\[
x\in[-1,1],
\qquad
\theta\in\R/2\pi\Z,
\qquad
t\in\R/L\Z.
\]
The metric
\begin{equation}\label{eq:intro-gstar}
g_*=\frac{5\,dx^2}{(1-x^2)(4+x^2)}
+\frac{(1-x^2)(4+x^2)}5\,d\theta^2
+\frac{(4+x^2)^2}{16}\,dt^2
\end{equation}
is smooth, has scalar curvature $R\ge\frac25$, and satisfies
\[A(\mathbb S^2\times \mathbb S^1)>-1,\]
and \begin{equation}
 \int_M |\Ricz|^2\dv<\tfrac16\int_M R^2\dv.
\end{equation}
\end{example}

\begin{proof}
 We use \eqref{eq:moment} with
\begin{equation}\label{eq:counter-first-metric}
 F=\tfrac15(1-x^2)(4+x^2),\qquad b=1+\tfrac14x^2.
\end{equation}
The function $F$ satisfies \eqref{eq:c2}.
Hence the metric $g_*$ in \eqref{eq:intro-gstar} extends smoothly over the
poles. From \eqref{eq:counter-moment-curvature} we obtain
\begin{equation}\label{eq:counter-first-curvature}
 \begin{aligned}
 K&=\tfrac35(1+2x^2),\\
 \Delta b&=\tfrac25-\tfrac9{10}x^2-\tfrac12x^4,\\
 \frac{|\hesz b|^2}{b}&=\frac{\frac{1}{2} F^2\left(b^{\prime \prime}\right)^2}{b}=\tfrac1{50}(1-x^2)^2(4+x^2),\\
 \langle\nabla_h K,\nabla_h b\rangle_h&=FK'b'=\tfrac65Fx^2.
 \end{aligned}
\end{equation}
Inserting \eqref{eq:counter-first-curvature} into the scalar curvature formula in \eqref{eq:warped4}, we obtain
\begin{equation}\label{eq:counter-first-scalar}
 R=\frac{2(16x^4+45x^2+4)}{5(x^2+4)} \ge \frac 25.
\end{equation}

For the metric $g_*$, we have $dv_{g_*}=b\,dx\,d\theta\,dt$.
One can compute 
\begin{equation}
 \begin{aligned}
\int_{-1}^{1}b\,dx&=\tfrac{13}{6},&
\int_{-1}^{1}Kb\,dx&=\tfrac{111}{50},&
 \bar{K}_b=\frac{333}{325},
 \end{aligned}
\end{equation}
and
\begin{align*}
2\int_{-1}^1 FK'b'\,dx&=\frac{496}{875},\\
\int_{-1}^1 b(K-\bar K_b)^2\,dx&=\frac{16392}{56875},\\
\int_{-1}^1\frac{F^2(b'')^2}{b}\,dx&=\frac{464}{2625}.
 \end{align*}
Thus
\[
A(M,g)+1
=\frac{13/6}{(111/50)^2}
\left(\frac{496}{875}-\frac{16392}{56875}-\frac{464}{2625}\right)
=\frac{34768}{776223}>0.
\]
Thus $A(\mathbb S^2\times \mathbb S^1,g_*)=-741455/776223>-1$ on $\mathbb S^2\times \mathbb S^1$, and thus, $A(\mathbb S^2\times \mathbb S^1)>-1$.

By
\[
|\Ricz|^2
=\frac{|\hesz b|^2}{b^2}
 +\frac16\left(2K+\frac{\Delta_h b}{b}\right)^2.
\]
Substituting $F$ and $b$, we obtain
\begin{align*}
b|\Ricz|^2
&=\frac{F^2(b'')^2}{2b}
  +\frac b6\left(2K+\frac{\Delta_h b}{b}\right)^2\\
&=\frac{2x^6+19x^4}{75}
  +\frac{69x^2}{50}-\frac{54}{25}
  +\frac{32}{3(x^2+4)},
\end{align*}
and hence
\[
\frac{\displaystyle\int_M |\Ricz|^2\,dv_{g_*}}
     {\displaystyle\int_M R^2\,dv_{g_*}}
=\frac{\int_{-1}^1 b|\Ricz|^2\,dx}
       {\int_{-1}^1 bR^2 \,dx}
=\frac{8639-28000\arctan\frac12}{6\,(45963-112000\arctan\frac12)}\approx0.1213<\frac16.
\]

\end{proof}

The functions $F$ and $b$ in \eqref{eq:counter-first-metric} are even, and the antipodal
map $(x,\theta)\mapsto(-x,\theta+\pi)$ is an isometry of $h$.
Writing $p\mapsto-p$ for the antipodal map, we obtain three isometric involutions of $g_*$,
\begin{equation}\label{eq:counter-involutions}
 \begin{aligned}
 (p,t)&\mapsto(-p,t),\\
 (p,t)&\mapsto(-p,-t),\\
 (x,\theta,t)&\mapsto(x,-\theta,t+L/2).
 \end{aligned}
\end{equation}
The corresponding
quotients are $\mathbb{RP}^2\times \mathbb S^1$, $\mathbb{RP}^3\#\mathbb{RP}^3$,
and the non-orientable two-sphere bundle over the circle, on which $A(M)>-1$ since the quotient
$A(g)$ is invariant under finite Riemannian covers. By \cite{To},
 $\mathbb S^2\times \mathbb S^1$, $\mathbb{RP}^2\times \mathbb S^1$, $\mathbb{RP}^3\#\mathbb{RP}^3$,
and the non-orientable two-sphere bundle over the circle are the only four manifolds diffeomorphic to the closed quotients of $\mathbb S^2\times \mathbb R$.

\

\noindent {\it AI disclosure.} During the preparation of this article, the authors used OpenAi Astra 6  to explore potential proof strategies and construct the counterexamples.  All mathematical statements and final proofs were independently verified and rigorously proved by the authors, who take responsibility for the correctness of the paper.

\bibliography{integralpinching}
\bibliographystyle{amsplain}
\end{document}